\documentclass[12pt,a4paper]{article}
\usepackage{graphicx}
\usepackage{amsfonts}
\usepackage{float}
\usepackage[latin1]{inputenc}
\usepackage{blkarray}
\usepackage{color}
\usepackage{amsmath}
\usepackage{amstext}
\usepackage{amssymb}
\usepackage{array}%,arydshln}
\usepackage{multicol}
\usepackage{pdfpages}
\usepackage{multirow}
\usepackage{threeparttable}
\usepackage{url}
\usepackage{enumerate}
\usepackage{xcolor}
\usepackage[ruled,vlined]{algorithm2e}
\usepackage{algorithmic}

\usepackage[bbgreekl]{mathbbol}
\DeclareSymbolFontAlphabet{\mathbb}{AMSb}
\DeclareSymbolFontAlphabet{\mathbbl}{bbold}

\usepackage{subcaption}
\def\R{\Bbb R}

\def\<{\langle}
\def\>{\rangle}

\def\wt{\widetilde}
\def\Chi{\raise .3ex \hbox{\large $\chi$}}

\def\mass{{\rm h_2o} }
\def\energy{e}
\def\WI{W\!\!I}

\def\[{\Bigl [}
\def\]{\Bigr ]}
\def\({\Bigl (}
\def\){\Bigr )}
\def\[{\Bigl [}
\def\]{\Bigr ]}
\def\({\Bigl (}
\def\){\Bigr )}
\def\l{\iota}

\def\div{{\mbox{\rm div}}}
\def\dsp{\displaystyle}
\def\K{\mathbb K}

\def\s{{\bf v}}

\def\cells{{\cal M}}

\def\nodes{{\cal V}}
\def\edges{{\cal E}}

\def\welledge {{\mathfrak a}}

\def\l{{\ell}}
\def\g{{\rm g}}

\begin{document}

%%-----------------------------
%%      the top matter
%%-----------------------------
\title{Multi-segmented non-isothermal compositional liquid gas well model for geothermal processes}

\author{ 
D. Castanon Quiroz \thanks{Instituto de Investigaciones en Matem\'aticas Aplicadas y en Sistemas, Universidad Nacional Aut\'onoma de M\'exico, Circuito Escolar s/n, Ciudad Universitaria C.P. 04510 Cd. Mx., M\'exico, daniel.castanon@iimas.unam.mx},
L. Jeannin \thanks{STORENGY, 12 rue Raoul Nordling - Djinn - CS 70001 92274 Bois Colombes Cedex, France, laurent.jeannin@storengy.com}, 
S. Lopez\thanks{BRGM, 3 avenue Claude-Guillemin, BP 36009, 45060 Orl\'eans Cedex 2, France, s.lopez@brgm.fr}, 
R. Masson\thanks{Universit\'e C\^ote d'Azur, Inria, CNRS, LJAD, 
UMR 7351 CNRS, team Coffee, Parc Valrose 06108 Nice Cedex 02, France, roland.masson@univ-cotedazur.fr} 
}

\maketitle

\abstract{
  We consider a non-isothermal compositional gas liquid  model for the simulation of well operations in geothermal processes. The model accounts for phase transitions assumed to be at thermodynamical equilibrium and is based on an hydrodynamical Drift Flux Model (DFM) combined with a No Pressure Wave approximation of the momentum equation.
  The focus of this work is on the design of a robust discretization accounting for slanted and multibranch wells with the ability to simulate both transient behavior such as well opening as well as coupled simulations at the time scale of the reservoir.
  It is based on a staggered finite volume scheme in space combined with a fully implicit Euler time integration. 
The construction of consistent and stable numerical fluxes is a key feature for a robust numerical method. It is achieved by combining 
a monotone flux approximation for the phase superficial velocities  with an upwind approximation of the phase molar fractions, density and enthalpy. In order to facilitate the coupling of the well and reservoir models, the Newton linearization accounts for the elimination of the hydrodynamical unknowns leading to Jacobian systems using the same primary unknowns than those of the reservoir model.
The efficiency of our  approach is investigated on both stand alone well test cases without and with cross flow, and on a fully coupled well-reservoir simulation.  
}


%%-----------------------------
%%      your text
%%-----------------------------

\section{Introduction}

Wells are key objects during the operation of geothermal reservoirs: they provide the connection between the surface and the geothermal reservoir at depth, enable fluid and heat exchanges with the geological structure, and are the place of two-phase flows and phase changes. During apraisal phases well tests provide valuable information on the target and during geothermal field operation, it is very important to correctly describe the behavior of  wellbore flows in order to predict  the quality of the fluids produced at the surface. A thorough understanding of wells and their modeling allow a better understanding of the reservoir behaviour.

Physical phenomena and exchanges with the reservoir are complex. For example, the start-up of a high-temperature geothermal well, fed by a hot liquid, gives rise to transient flows with flash in the well. Geothermal wells may also be connected by several feed zones at different pressures and some feed zones will produce in the well, while others may receive fluids from the well, even when the latter is in production.  Finally, in addition to describing wells over transient short periods, well modeling can also be used to predict the operational behaviour of a reservoir over a period of several years. In this case, the well model has to be coupled with a flow model in the geological structures.

Flows occuring in geothermal wells are generally described as  two-phase thermal pipe flows.  Given the well's radius/length aspect ratio, the usual assumption is to rely on 1D models. This type of flow is relevant not only to the geothermal industry, but also to the oil and gas industry (for example, \cite{shippen_steady-state_2012} for steady flows, \cite{danielson_transient_2012} for transient flows or \cite{jerez-carrizales_prediction_2015}). 
The nuclear industry has also a strong interest in multiphase thermal flows in pipes (\cite{delhaye_thermohydraulique_2008} for a review).
 We present below a few bibliographical elements relevant to our work for geothermal processes, but readers interested in more details may find further information in \cite{tonkin_review_2021}, where Tonkin et al. recently proposed an in-depth review of mathematical models of flows in geothermal boreholes.

The specifics of geothermal flows lie in the description of the interaction between the well and the reservoir, and the consideration of feed zones and cross-flow phenomena. Bjornsson \cite{bjornsson_multi-feedzone_1987} was one of the first to propose a multi-feedzone well simulator. %using HOLA software%. 
These well models are either transient or steady-state; transient well models are particularly useful for interpreting well tests or, for example, well start-up, while  steady-state models are more dedicated to operations follow-up  (for example, \cite{gudmundsdottir_wellbore_2015} for matching field data or \cite{battistelli_modeling_2020} for the modeling of the IDDP-1 supercritical well).

One of the greatest difficulties in modeling multiphase flows in geothermal wells arises from the complexity and variety of flow regimes encountered. The simplest homogeneous flow model, in which all phases flow at the same velocity, is unsuitable because it cannot satisfactorily reproduce the in-situ well volume fraction and flow rate of each phase. Drift Flux Models (DFM) are more complex models in which a drift velocity (or slip law) is introduced to describe the relative motion between phases  \cite{zuber_average_1965, wallis_one-dimensional_1969}. These DFM require a number of empirical correlations derived from experiments for different flow regimes in order to correctly describe phase slip and in situ phase volume fractions  \cite{wallis_one-dimensional_1969, delhaye_thermohydraulique_2008, tonkin_review_2021}.
Given the diversity of flow regimes encountered in geothermal wells, preference is given to DFM capable of describing all flow regimes encountered during production well operations \cite{tonkin_review_2021, shi2005}.

It is usually assumed that the thermodynamical equilibrium assumption applies in the well, but some authors  are interested in non-equilibrium phenomena \cite{pereira_exsolving_2020, aunzo_wellbore_1991} to describe scaling or kinetics phenomena.
\cite{satman_effect_1999, lopez_40_2010, khasani_numerical_2021, zolfagharroshan_rigorous_2020}. Nevertheless, these reactive phenomena will not be considered in our work. Our objective is to set-up an extensible framework for the modeling of transient multi-component two-phase geothermal wellbore flows with a focus on robust numerical schemes and the ability to couple reservoir and wellbore flow models.

The considered well thermodynamical model accounts for a two-phase liquid gas compositional non-isothermal system as for the reservoir model described in \cite{Xing.ea:2017}. This two-phase compositional system uses a Coats-type formulation in order to account for a large class of equilibrium thermodynamical models ranging from single phase, two-phase immiscible to two-phase fully compositional with phase appearance and disappearance. Choosing the same formulation for both the well and the reservoir ensures robust coupling between these two systems. Then, following \cite{Xing2018}, to be able to consider complex well geometries with slanted or multi-branched paths, the geometry of each well is defined as a set of edges of the reservoir mesh assumed to define a rooted tree oriented away from the root.

The hydrodynamical model is based on the DFM introduced in \cite{shi2005} for the slip law. Moreover, following \cite{Faille02,LIVESCU2010138,masella_transient_1998}, we consider a" No Pressure Wave" (NPW) approximation, where the mixture momentum equation is reduced to a static balance accounting for gravity and wall friction pressure losses. This assumption remains valid to describe pressure transients (see \cite{tonkin_review_2021} for a numerical validation).

The discretization is based on a staggered finite volume scheme with  edge center superficial velocities and nodal pressure, temperature, saturation and molar fraction unknowns. This choice is consistent with the Vertex Approximate Gradient (VAG) nodal discretization \cite{Xing.ea:2017} used in the reservoir. It is combined with a two-point monotone flux for the phase superficial velocities and a phase based upwind approximation of the phase molar fractions, molar density and molar enthalpy. The monotonicity of the superficial velocity fluxes is a key property for the stability of the scheme. Moreover it guarantees the consistent definition of the phase based upwinding according to the sign of the phase superficial velocity.
The clear understanding of this interplay between the upwinding of the phase superficial velocities and of the upwinding of the phase transported variables is one of the main contribution of this work. 
The monotonicity property of the superficial velocites is obtained on the full range of saturations $[0,1]$ using an hybrid upwinding approach exploiting the mathematical properties of the drift velocity and profile parameter of the DFM model proposed in \cite{shi2005}. Note that the methodology developped in this work could be applied to other slip laws sharing the same types of mathematical properties. 

The time integration is fully implicit to account for the large reservoir time scale and for the stiffness of the system at the gas phase appearance. 
The Jacobian system is reduced by elimination of the phase superficial velocities and of the well flow rates using the slip law, the momentum NPW equation and the well monitoring conditions. This elimination allows the use of the same primary unknowns as for the reservoir model  \cite{Xing.ea:2017} which facilitates the fully implicit coupling with the reservoir system.

The remainder of this paper is organized as follows. Section \ref{sec:wellmodel} presents the well physical model and its staggered finite volume discretization is described in Section \ref{subsec_well}. It starts, in Subsection \ref{discrete_mesh_not}, by the definition of the well mesh and discrete unknowns, while the staggered finite volume scheme is detailed in Subsection \ref{discrete_conseq}. Subsection \ref{discrete_velocities} defines the monotonicity properties required for the phase superficial velocities and underlines the interplay between these properties and the upwind approximation of the phase molar fractions, density and enthalpy. The detailed specific construction of the monotone superficial velocity flux function based on an hybrid upwinding approach for the DFM model proposed in \cite{shi2005} is reported to Appendix \ref{sec_monotoneflux}. The well monitoring conditions are introduced in Subsection \ref{sec:monitoring} in the case of a production well and Appendix \ref{sec:peaceman} recall the definition of the reservoir-well fluxes based on the Peaceman indexes.  
Then, the algorithm used to solve the coupled nonlinear system at each time step of the simulation is addressed in Section \ref{sec:nonlinear}. It is based on a Newton-Raphson algorithm combined with an active set method for the phase appearance/disappearance and the well monitoring conditions. Its main originality is related to the elimination in the linearization process of the thermodynamical, hydrodynamical and monitoring closure laws leading to the same set of primary unknowns as for the reservoir model in the Jacobian system. This procedure is described respectively in  Subsections \ref{sec:elimthermo}, \ref{sec:elimhydro} and \ref{sec:elimmonitoring}. 
The objective of the numerical Section \ref{sec:testcases} is to validate the discrete model and to investigate the efficiency of our approach on both stand alone well test cases in Subsection \ref{sec:standalone} and and for a fully coupled well-reservoir model in Subsection \ref{sec:fullycoupled}.  The first test case in Subsection \ref{sec:buckleyleverett} considers the validation of the numerical model on a simple two-phase incompressible immiscible flow which reduces to a Buckley Leverett scalar hyperbolic equation.
The second and third test cases consider a high energy liquid vapor single component thermal flow along a multi-branch production well without and with cross flow in respectively Subsections \ref{sec:chairwell.test} and \ref{sec:crossflow}.
Finally, a single well liquid vapor geothermal test case is considered in Subsection \ref{sec:fullycoupled} to validate the fully coupled model by comparison to a simpler well model based on a single implicit well unknown.

\section{Well physical model}\label{sec:wellmodel}

Let $\mathcal{W}$ denote the set of wells. 
As in \cite{Xing2018}, each multi-branch well $\omega\in \mathcal{W}$ is defined by a set of oriented edges of the reservoir mesh assumed to define a rooted tree oriented away from the root. This orientation corresponds to the drilling direction of the well.
We consider a two-phase liquid gas, compositional, and non-isothermal flow model. 
The liquid ($\l$) and gas ($\g$) phases are described by their pressure $p$ (both liquid and gas pressures are assumed to match along the well), temperature $T$, volume fraction or saturation $s^\alpha$, and molar fractions $c^\alpha=(c_i^\alpha)_{i\in \mathcal{C}}$, $\alpha\in \mathcal{P} = \{\l,\g\}$, where $\mathcal{C}$ denotes the set of components and $\mathcal{P}$ the set of liquid and gas phases. The components are not necessarily present in both phases but at least in one of them. We consequently denote by $\mathcal{C}^\alpha\subset \mathcal{C}$ the set of components in phase $\alpha$ and by $\mathcal{P}_i$ the set of phases containing the component $i\in \mathcal{C}$.

The thermodynamical model uses a Coats' formulation already developped for the reservoir model (see e.g. \cite{Xing.ea:2017}). It is based on the natural variables $p$, $T$, $s^\alpha$, $c^{\alpha}$, $\alpha \in Q$ with a set of unknowns and equations depending on the additional unknown $Q\subset \mathcal{P}$ representing the subset of present phases at each point of the space time domain. This formulation has the advantage to account for an arbitrary number of phases and components and also to allow the components to be either present or absent of each given phase. This allows to account in the same framework for models ranging from single phase liquid, single phase gas to two-phase gas liquid, and from immiscible to fully compositional. 

For each phase $\alpha$, we denote by
 $\zeta^\alpha(p,T,c^\alpha)$ its molar density, by 
 $\rho^\alpha(p,T,c^\alpha)$ its specific density, by 
 $\mu^\alpha(p,T,c^\alpha)$  its dynamic viscosity, by 
 $e^\alpha(p,T,c^\alpha)$ its molar internal energy, and by 
$h^\alpha(p,T,c^\alpha)$ its molar enthalpy.
The fugacities are  denoted by $f_i^\alpha(p,T,c^\alpha)$. \\

The hydrodynamical model is based on the Drift Flux Model (DFM) with slip closure laws expressing the phase superficial velocities as functions of the mixture velocity and drift velocity terms \cite{ZF65}. It is combined with the No-Pressure-Wave (NPW) approximation of the momentum equation relating the pressure gradient to the friction and gravity pressure losses \cite{Faille02,LIVESCU2010138}.
 
Let $\tau$ denote the curvilinear coordinate along the well and ${\bf e}_\tau$ the unit tangential vector along the well assumed to be oriented toward the well root node. We denote by ${\bf u}^\alpha = u^\alpha {\bf e}_\tau$ the superficial velocity of each phase $\alpha\in \mathcal{P}$.
The mixture velocity is defined by
$$
{\bf u}^m = {\bf u}^\l + {\bf u}^\g = u^m {\bf e}_\tau. 
$$
The tangential divergence along the well is denoted by $\div_\tau$ and the tangential gradient by $\nabla_\tau$. The section along the well is denoted by $S^\omega(\tau)$. We also set ${\bf g}_\tau = ({\bf g}\cdot {\bf e}_\tau){\bf e}_\tau $ where ${\bf g}$ is the acceleration of gravity vector.\\

Let $\epsilon_\omega(\tau) = {(\nabla_\tau z(\tau))\cdot {\bf e}_\tau) \over |\nabla_\tau z(\tau)|} \in \{-1,1\}$ define the orientation along the well, with arbitrary value for $\nabla_\tau z=0$ along an horizontal branch. Following  \cite{shi2005}, let us introduce the function 
$$
o_\omega(\tau) = \epsilon_\omega(\tau) \cos(\theta_\omega(\tau))^{1/2} (1+ \sin(\theta_\omega(\tau)))^2,
$$
where $\theta_\omega(\tau) \in [0,\pi /2]$ is the edge angle w.r.t. the vertical direction. Note that $o_\omega(\tau) = 0$ along an horizontal branch. 
The  slip closure law introduced in \cite{ZF65} expresses the gas superficial velocity $u^\g$  as a function of the mixture velocity, the drift velocity, and the profile parameter: 
\begin{equation}
  \label{slip_law}
u^\g  =   s^g U_{\rm d}(s^\g,\bar H) o_\omega +   s^g C_0(s^\g,\bar H) u^{\rm m},
\end{equation}
with the non negative drift velocity $U_{\rm d}$ (for vertical wells) and the profile parameter $C_0$ such that $C_0(s^\g=1,\bar H) = 1$, $U_{\rm d}(s^\g =1,\bar H) = 0$, $s^\g C_0(s^\g,\bar H) \leq 1$ and
$s^\g C_0(s^\g,\bar H)$ non decreasing w.r.t. $s^\g$.
The drift velocity and profile parameter depend on a set of thermodynamical variables denoted by
$\bar H$ and typically comprising the densities of the phases and the liquid gas interfacial tension. They also depend on other fixed parameters such as the well radius $r_\omega$. \\

We denote by  
$$
{\bf T}^f = T^f {\bf e}_\tau,
$$
the wall friction law depending on the mixture velocity ${\bf u}^m$ and on thermodynamical quantities. To fix ideas, we will use in the following a Darcy-Forchheimer type law:
\begin{equation}
  \label{eq_darcy_forchheimer}
T^f = - \( { 8 \mu^{\rm m} \over r^2_\omega} + f_q \rho^{\rm m} { 1\over 4 r_\omega} |u^{\rm m} | \) u^{\rm m},  
\end{equation}
where $f_q$ is the friction coefficient, $\mu^{\rm m} = s^\g \mu^\g + s^\l\mu^\l$  the mean viscosity and $\rho^{\rm m} = s^\g \rho^\g + s^\l\rho^\l$  the mean specific density. \\

The well model is based on the molar conservation  of each component \eqref{eq_cons_mole_well}, the total energy conservation \eqref{eq_cons_energy_well} and  the mixture momentum equation based on the NPW model \eqref{eq_momentum_well}.
It is combined with the the DFM slip law \eqref{eq_slip_well} and the thermodynamical equilibrium \eqref{eq_thermo_well}:     
\begin{subequations}\label{well_model}
\begin{alignat}{2}
  & \dsp \sum_{\alpha\in \mathcal{P}_i} \(\partial_t(S^\omega c_i^{\alpha} \zeta^\alpha s^{\alpha}) + \div_\tau ( S^\omega  c_i^{\alpha} \zeta^\alpha {\bf u}^\alpha) =  q^{r\rightarrow \omega}_{i}, \, i\in{\cal C}, \label{eq_cons_mole_well}\\
  & \dsp \sum_{\alpha\in \mathcal{P}} \( \partial_t ( S^\omega  e^\alpha \zeta^\alpha s^{\alpha}) + \div_\tau ( S^\omega  h^\alpha \zeta^\alpha {\bf u}^\alpha) \) +\div_\tau( -S^\omega  \lambda \nabla_\tau T) =  q^{r\rightarrow \omega}_{e},
  \label{eq_cons_energy_well} \\
& \dsp \nabla_\tau p = {\bf T}^f  + \rho^{\rm m} ~{\bf g}_\tau, \label{eq_momentum_well}
 \end{alignat}
\begin{equation}\label{eq_thermo_well}
\begin{aligned}  
  & f^g_i(p,T,c^{g}) = f^l_i(p,T,c^{l}) \, \mbox{ if } Q\cap \mathcal{P}_i =\mathcal{P},\\
  & \dsp \sum_{i\in {\cal C}^\alpha} c_i^{\alpha} = 1,  \alpha\in Q,\\
  & \dsp \sum_{\alpha \in Q} s^\alpha = 1, \quad s^{\alpha} = 0 \, \mbox{ if } \alpha\not\in Q,  
\end{aligned}
\end{equation}
\begin{equation}  \label{eq_slip_well}
\dsp u^\g  =   s^g U_{\rm d}(s^\g,\bar H) o_\omega +   s^g C_0(s^\g,\bar H) u^{\rm m} \,    \mbox{ if } Q=\mathcal{P},
\end{equation}
\end{subequations} 
where $q^{r\rightarrow \omega}_{i}$ and $q^{r\rightarrow \omega}_{e}$ are the molar and energy exchange terms with the reservoir. 
The system is closed with monitoring conditions at the root node and with a flash computation to determine the set of present phases $Q$. The temperature and pressure continuity is also assumed at the well junctions.

\section{Staggered Finite Volume discretization of the well model}
\label{subsec_well}

Notations about the well mesh and the discrete variables are introduced in Subsection \ref{discrete_mesh_not}. 
The scheme presented in Subsection \ref{discrete_conseq} is based on a fully implicit time integration to cope with large time steps at the reservoir time scale. It is combined with a staggered finite volume discretization in space using node centred control volumes for the molar and energy conservation equations and edge control volumes for the momentum equation. A key ingredient is the discretization of the convective fluxes based on a monotone two-point flux for the approximation of the superficial velocities. This framework is presented in Subsection \ref{discrete_velocities} and an example is detailed in Appendix \ref{sec_monotoneflux} based on the model proposed in \cite{shi2005}. Thanks to the monotonicity and consistency properties of the superficial velocities, it is combined with an upwind approximation of the phase molar fractions, density and enthalpy w.r.t. the sign of the phase superficial velocity. The monitoring conditions at the head node of the well are described in Subsection \ref{sec:monitoring} in the case of a production well considered in the numerical Section.

\subsection{Well discretization and notations} \label{discrete_mesh_not}

\begin{figure}[H]
\begin{center}
\includegraphics[width=0.35\textwidth]{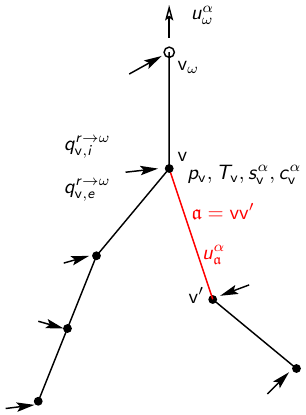}   
\caption{Example of multi-branch well $\omega$ with its root node $\s_\omega$, one edge $\mathfrak{a}=\s\s'$ ($\s$ parent node of $\s'$) and the main physical quantities: the head node superficial velocities $u^\alpha_\omega$ (non negative for production wells and non positive for injection wells), the molar and energy flow rates between the reservoir and the well $q^{r\rightarrow \omega}_{\s,i}$, $q^{r\rightarrow \omega}_{\s,\energy}$, the well node pressure, temperature, saturations and phase molar fractions  $p_\s, T_\s, s^\alpha_{\s},c^{\alpha}_{\s}$, and the phase superficial velocity $u_\welledge^\alpha$ at the edge $\welledge$ oriented positively from ${\s'}$ to ${\s}$. }
\label{fig_wellmodel}
\end{center}
\end{figure}

  The set of nodes of a well $\omega\in \mathcal{W}$ is denoted by $\nodes_\omega$ and
its root node is denoted by $\s_\omega$.  A partial ordering is defined on the set of vertices $\mathcal{V}_\omega$
with $\s \underset{\omega}{<} \s'$ if and only if the unique path from the root $\s_\omega$ to $\s'$ passes through $\s$.
The set of edges of the well $\omega$ is denoted by $\mathcal{E}_\omega$ and for each edge $\welledge \in \edges_\omega$
we set $\welledge={\s}{\s'}$ with ${\s} \underset{\omega}{<} {\s'}$ (i.e. $\s$ is the parent node of $\s'$, see Figure \ref{fig_wellmodel}). It is assumed that $\nodes_{\omega_1}\cap \nodes_{\omega_2} = \emptyset$ for any $\omega_1, \omega_2 \in \mathcal{W}$ such that $\omega_1 \neq \omega_2$.

Let $|S^\omega_\welledge|$ and $|S^\omega_{\s_\omega}|$ denote the well section at respectively the center of the edge $\welledge$  and at the head node $\s_\omega$. For $\welledge\in \mathcal{E}_\omega$, $|\welledge|$ denotes the lengh of the edge $\welledge$. Let us define $\edges_\s^\omega\subset \edges_\omega$ as the set of well edges sharing the node $\s\in \nodes_\omega$. For all ${\s}{\s'} = \welledge\in \edges_\omega$, let us set $\kappa_{\welledge,\s'} = -1$ and $\kappa_{\welledge,\s} = 1$.
For each edge $\welledge = {\s}{\s'}\in \mathcal{E}_\omega$, let us denote by $u^\alpha_\welledge$ the superficial velocity of phase $\alpha$ along the edge $\welledge$ oriented positively from ${\s'}$ to ${\s}$. Using the Coats' formulation, the set of well unknowns is defined 
at each node $\s\in \nodes_\omega$ by
\begin{itemize}
\item the set of present phases $Q_\s$,  
\item the well pressure $p_\s$,
\item  the well temperature $T_\s$,
\item the well saturations $s_{\s}^\alpha$ for  $\alpha\in \mathcal{P}$, 
\item the well molar fractions $c^{\alpha}_\s$ for  $\alpha\in Q_\s$,
\item the number of moles $\tilde n_{\s,i}$ for $i\in \wt{\mathcal{C}}_{Q_\s}$,
\end{itemize}
and by 
\begin{itemize}
\item the edge superficial velocities $u^\alpha_\welledge$ for $\alpha \in \mathcal{P}$ at each well edge $\welledge \in \edges_\omega$,
\item the head node superficial velocities  $u^\alpha_\omega$ for $\alpha \in \mathcal{P}$ at the well head node $\s_\omega$.  
\end{itemize}
In the above definition, $\wt{\mathcal{C}}_{Q}$ denotes the set of components not contained in $\bigcup_{\alpha\in Q} \mathcal{C}^\alpha$. For $\wt{\mathcal{C}}_{Q_\s}\neq \emptyset$, the additional unknowns $\tilde n_{\s,i}$ corresponding to the number of moles of the absent components $i\in \wt{\mathcal{C}}_{Q_\s}$ are needed first to avoid the singularity of the Jacobian and second to track the appearance of the missing phase containing these absent components (see \cite{Xing.ea:2017}  for details). 
Note that the saturation $s_{\s}^\alpha$ of an absent phase $\alpha \not \in Q_\s$ vanishes. 
In the following, to fix ideas, we assume that the monitoring conditions at the well head node correspond to the case of a production well consistently with the numerical section. \\

For a given thermodynamical phase property $\xi^\alpha$  at a node $\s$, we will use the notation
$$
\xi^\alpha_\s = \xi^\alpha(p_\s,T_\s,c^\alpha_\s), 
$$
with typically $\xi^\alpha = \zeta^\alpha, \rho^\alpha, e^\alpha, h^\alpha, f^\alpha_i, \mu^\alpha$.  \\

The molar and energy flow rates between the reservoir and the well $\omega$ at a given node $\s \in \mathcal{V}_\omega$ are defined by a two point flux approximation between the reservoir and well properties at node $\s$ based on the Peaceman approach combined with a phase potential upwinding. They are denoted respectively by $q^{r\rightarrow \omega}_{\s,i}$ for each component $i\in \mathcal{C}$ and $q^{r\rightarrow \omega}_{\s,e}$ in what follows (see Figure \ref{fig_wellmodel} and Appendix \ref{sec:peaceman} for their definitions).

\subsection{Fully implicit Finite Volume scheme}
\label{discrete_conseq}

Let us define at each nodal control volume $\s$, the total number of mole of each component and the total energy as follows 
$$
n_{\s,i} =
\left\{\begin{array}{ll}
\dsp ( \sum_{\welledge\in \edges_\s^\omega} |S^\omega_\welledge| {|\welledge| \over 2}) \sum_{\alpha\in Q_\s\cap \mathcal{P}_i} c_{\s,i}^{\alpha} s^\alpha_{\s} \zeta^\alpha_\s & \mbox{ if } i \not\in \wt{\mathcal{C}}_{Q_\s},\\[4ex]
\dsp \tilde n_{\s,i} & \mbox{ else},
\end{array}\right.
$$
$$
n_{\s,e} = (\sum_{\welledge\in \edges_\s^\omega} |S^\omega_\welledge| {|\welledge| \over 2}) \sum_{\alpha\in Q_\s} s^\alpha_{\s} \zeta^\alpha_\s  e^\alpha_\s, 
$$
where the case of an absent component $i\in \wt{\mathcal{C}}_{Q_\s}$ is accounted for by the introduction of the additional unknown $\tilde n_{\s,i}$ in order to avoid the singularity of the discrete system. 

For each edge $\welledge={\s}{\s'}\in \edges_\omega$, and each phase $\alpha$, let us define the following upwind approximation of the phase molar fractions, density and enthalpy w.r.t. the sign of the phase superficial velocity: 
\begin{eqnarray}
\label{eq_ha_rhoa_sa}
 \xi^\alpha_\welledge =
\left\{\begin{array}{r@{\,\,}c@{\,\,}l}
& \xi^\alpha_{\s'} \mbox{ if } u^\alpha_\welledge \geq 0,\\[1ex]
& \xi^\alpha_{\s} \mbox{ if } u^\alpha_\welledge < 0,  
\end{array}\right.
\end{eqnarray}
with $\xi^\alpha = c^\alpha, \zeta^\alpha, h^\alpha$. \\

We denote by $(0,t_F)$ the time interval and consider its discretization given by $t^n$, $n=0,\cdots,N$ with $t^0 = 0$, $t^N = t_F$ and $\Delta t^n = t^n - t^{n-1} > 0$ for all $n=1,\cdots,N$. The time discretization of 
the system \eqref{well_model} is based on an implicit Euler integration scheme. To simplify the notations only the accumulation terms $n_{\s,i}, n_{\s,e}$ at time $t^{n-1}$  will be specified by the $n-1$ superscript. For all other quantities considered at the current time $t^n$, we will drop the $n$ superscript for simplicity. 

Then, the discretization of the well equations \eqref{eq_cons_mole_well}-\eqref{eq_cons_energy_well}-\eqref{eq_momentum_well}-\eqref{eq_thermo_well} at each time step $n=1,\cdots,N$ is defined as follows: 
\begin{subequations}
  \label{eq_well_discrete}
  \begin{equation} \label{cons_molar_discrete}
    \begin{aligned}
      \dsp {n_{\s,i}- n_{\s,i}^{n-1}\over \Delta t^n }  & + \sum_{\welledge\in \edges_\s^\omega}  \sum_{\alpha\in \mathcal{P}_i} - \kappa_{\welledge,\s} |S^\omega_\welledge| c_{\welledge,i}^\alpha  \zeta^\alpha_\welledge u^\alpha_\welledge \\
      & = \dsp q^{r\rightarrow \omega}_{\s,i} - \delta_\s^{\s_\omega}  \sum_{\alpha\in \mathcal{P}_i}  |S^\omega_{\s_\omega}| c^{\alpha}_{\s,i} \zeta^\alpha_\s u^\alpha_\omega,  \, i\in \mathcal{C}, \qquad \s \in \mathcal{V}_\omega,
       \end{aligned}
 \end{equation}
 \begin{equation}\label{cons_energy_discrete}
   \begin{aligned}
   \dsp {n_{\s,e}- n_{\s,e}^{n-1}\over \Delta t^n }  & + \sum_{\welledge\in \edges_\s^\omega} |S^\omega_\welledge| \( { \lambda_\welledge \over |\welledge|} (T_\s - T_{\s'}) + \sum_{\alpha\in \mathcal{P}  } - \kappa_{\welledge,\s}  h^\alpha_\welledge \zeta^\alpha_\welledge u^\alpha_\welledge  \)  \\
     & \dsp  = \dsp q^{r\rightarrow \omega}_{\s,e} - \delta_\s^{\s_\omega} \sum_{\alpha\in \mathcal{P}} |S^\omega_{\s_\omega}| h^\alpha_\s \zeta^\alpha_\s u^\alpha_\omega, \qquad \s \in \mathcal{V}_\omega,
 \end{aligned}
 \end{equation}       
\begin{equation}\label{cons_momentum_discrete}
p_{{\s}} - p_{{\s'}} + \bar\rho^{\rm m}_\welledge g (z_{{\s}} - z_{{\s'}}) = T^f_\welledge |\welledge|, \qquad {\s}{\s'} = \welledge\in \edges_\omega,    
\end{equation}
\begin{equation}\label{eq_thermo_discrete}
\begin{aligned}  
& f^\g_{i,\s} = f_{i,\s}^\l  \mbox{ if } Q_\s\cap \mathcal{P}_i = \mathcal{P}, \qquad \s \in \mathcal{V}_\omega,\\   
& \dsp \sum_{i\in \mathcal{C}^\alpha} c^{\alpha}_{\s,i} = 1, \, \alpha\in Q_\s, \qquad \s \in \mathcal{V}_\omega,\\
& \dsp \sum_{\alpha\in Q_\s} s_{\s}^\alpha  = 1, \quad s^{\alpha}_\s = 0 \, \mbox{ if } \alpha\not\in Q_\s,  \qquad \s \in \mathcal{V}_\omega. 
\end{aligned}
\end{equation}
\end{subequations}
where $\delta$ stands for the Kronecker symbol with $\delta_\s^{\s_\omega} = 1$ if $\s=\s_\omega$, else $0$. 
In \eqref{eq_well_discrete}, the edge thermal conductivity of the mixture is defined by   
$$
\lambda_\welledge = \sum_{\alpha \in \mathcal{P}} {s^\alpha_{\s} + s^\alpha_{\s'} \over 2}  \lambda^\alpha, 
$$
with $\lambda^\alpha$ the thermal conductivity of the phase $\alpha$ assumed constant for simplicity. 
The wall friction term $T^f_\welledge$ is given by the Darcy-Forchheimer law \eqref{eq_darcy_forchheimer} leading to 
\begin{equation}
  \label{eq_darcy_forchheimer}
T^f_\welledge = - \( { 8 \bar\mu^{\rm m}_\welledge \over r^2_\welledge} + f_q \bar\rho^{\rm m}_\welledge { 1\over 4 r_\welledge} |u^{\rm m}_\welledge | \) u^{\rm m}_\welledge, 
\end{equation}
with the following arithmetic means of the mixture specific density and dynamic viscosity at the edge $\welledge$  
\begin{eqnarray}
  \label{eq_mum_rhom} 
\dsp \bar \rho^m_\welledge = {1\over 2} \( \sum_{\alpha\in Q_\s} s^\alpha_{\s} \rho^\alpha_\s + \sum_{\alpha\in Q_{\s'}} s^\alpha_{\s'} \rho^\alpha_{\s'} \), \qquad 
 \dsp \bar \mu^m_\welledge = {1\over 2} \( \sum_{\alpha\in Q_\s} \( s^\alpha_{\s} \mu^\alpha_\s + \sum_{\alpha\in Q_{\s'}} s^\alpha_{\s'} \mu^\alpha_{\s'} \). 
\end{eqnarray}

It remains to define the discretization of the edge superficial velocities $u^\alpha_\welledge$ accounting for the slip law \eqref{eq_slip_well} at given mixture velocity $u^{\rm m} = u^\g + u^\l$. This framework is detailed in the next subsection based on a monotone two-point flux. 

 \subsection{Monotone two-point flux for the edge superficial velocities} \label{discrete_velocities}

Let us define the arithmetic average of a thermodynamical variable $\xi^\alpha$ for a given phase $\alpha$ at a well edge $\welledge=\s\s'$ with $\alpha \in Q_\s \cup Q_{\s'}$ by 
$$
\bar \xi^\alpha_\welledge =
\left\{\begin{array}{ll}
\xi^\alpha(p_\s,T_\s,c^{\alpha}_{\s}) & \mbox{ if } \alpha \in Q_\s \mbox{ and } \alpha \not \in Q_{\s'},\\
\xi^\alpha(p_{\s'},T_{\s'},c^{\alpha}_{\s'}) & \mbox{ else if } \alpha \not\in Q_\s \mbox{ and } \alpha \in Q_{\s'},\\
\dsp { s^\alpha_{\s} \xi^\alpha(p_\s,T_\s,c^{\alpha}_{\s}) + s^\alpha_{\s'} \xi^\alpha(p_{\s'},T_{\s'},c^{\alpha}_{\s'}) \over s^\alpha_{\s} + s^\alpha_{\s'} } & \mbox{ else if }  s^\alpha_{\s} + s^\alpha_{\s'} > \epsilon,\\
\dsp {1\over 2} \(\xi^\alpha(p_\s,T_\s,c^{\alpha}_{\s}) + \xi^\alpha(p_{\s'},T_{\s'},c^{\alpha}_{\s'}) \) & \mbox{ else}. 
\end{array}\right.
$$
This arithmetic averaging is used to compute the thermodynamical variables $\bar H_\welledge$ at the well edge $\welledge=\s\s'$ which enter in the definition of the slip law \eqref{eq_slip_well}.

Then, the superficial velocities $u^\alpha_\welledge$ are obtained at given mixture velocity $u^{\rm m}_\welledge$ using a numerical two-point monotone flux denoted by
$$
F^\g_\welledge(s_{\s'}^\g,s_{\s}^\g,\bar H_\welledge),
$$
for the continuous flux function
$$
f_\welledge(s^\g) = s^\g U_{\rm d}(s^\g,\bar H_\welledge) o_\welledge +   s^\g C_0(s^\g,\bar H_\welledge) u^{\rm m}_\welledge, 
$$
where $o_\welledge$ is the value of the function $o_\omega$ along the edge $\welledge$. 
 The numerical flux function $F^\g_\welledge$ must satisfy the following properties:

\begin{itemize}
\item consistency property: 
$$
F_\welledge^\g(s^\g,s^\g,\bar H_\welledge) =  f_\welledge(s^\g) \mbox{ for all } s^\g\in [0,1], 
$$
\item monotonicity property which specifies that, for all $(s_{\s'}^\g,s_{\s}^\g) \in [0,1]\times [0,1]$,  
$F^\g_\welledge(s_{\s'}^\g,s_{\s}^\g,\bar H_\welledge)$ is non decreasing w.r.t. its first argument $s_{\s'}^\g$ and non increasing w.r.t. its second argument $s_{\s}^\g$.
\end{itemize}
Then the superficial velocities $u^\alpha_\welledge$ are defined by
\begin{equation}
  \label{eq_superficial_velocities}
\left\{\begin{array}{r@{\,\,}c@{\,\,}l}
u^\g_\welledge &=& F^\g_\welledge(s_{\s'}^\g,s_{\s}^\g,\bar H_\welledge),\\
u^\l_\welledge &=& u^{\rm m}_\welledge - u^\g_\welledge. 
\end{array}\right.
\end{equation}
From the consistency and monotonicity properties, it results that
\begin{itemize}
\item $s^\g_{\s'}=0 \Rightarrow u_\welledge^\g = F^\g(0,s^\g_{\s}) \leq 0$ for all $s^\g_\s \in [0,1]$, 
\item $s^\g_{\s}=0 \Rightarrow u_\welledge^\g = F^\g(s^\g_{\s'},0) \geq 0$ for all $s^\g_{\s'} \in [0,1]$, 
\item $s^\g_{\s'}=1 \Rightarrow u_\welledge^\l = u^{\rm m}_\welledge - F^\g(1,s^\g_{\s}) \leq 0$ for all $s^\g_\s \in [0,1]$, 
\item $s^\g_{\s}=1 \Rightarrow u_\welledge^\l = u^{\rm m}_\welledge - F^\g(s^\g_{\s'},1) \geq 0$ for all $s^\g_{\s'} \in [0,1]$. 
\end{itemize}  
These properties are key for the definition \eqref{eq_ha_rhoa_sa} of the upwind approximations of the phase molar fractions, phase molar density and phase molar enthalpy w.r.t. the sign of the phase superficial velocity. Indeed, they ensure that either the phase is present at the upwind node or that its superficial velocity vanishes ensuring that the product $\xi^\alpha_\welledge u^\alpha_\welledge$ is always properly defined, with $\xi^\alpha = \zeta^\alpha, c^\alpha, h^\alpha$.  \\

Note also that if $Q_\s \cup Q_{\s'} = \{\alpha\}$ (single phase $\alpha$ present at both nodes), the properties of the numerical flux function imply that $u_\welledge^\alpha = u^{\rm m}_\welledge$ and $u_\welledge^\beta = 0$ for the absent phase $\beta$. It means that, in that cases, the thermodynamical variables $\bar H_\welledge$ do not need to be computed which would have required a cumbersome extension of the molar fractions of the absent phase. \\

The specific construction of the two-point monotone flux function $F^g_\welledge$ depends on the laws $C_0$ and $U_{\rm d}$. We propose in Appendix \ref{sec_monotoneflux}  a monotone two-point flux for the gas liquid DFM model introduced in \cite{shi2005} on the full range of gas saturation.\\

\subsection{Monitoring conditions}\label{sec:monitoring}

Assuming the case of a production well to fix ideas, the monitoring conditions prescribe a minimum head node pressure $\bar p_\omega$ and a maximum well molar flow rate $\bar q_\omega \geq 0$.
Setting
\begin{equation}\label{eq_wellpwqw}
p_\omega = p_{\s_\omega} \mbox{ and } q_\omega = \sum_{\alpha\in \mathcal{P}} |S^\omega_{\s_\omega}| \zeta^\alpha_{\s_\omega} u^\alpha_\omega, 
\end{equation} 
the boundary condition at the head node $\s_\omega$ combines the complementary constraints on the pair $\(\bar q_\omega - q_\omega, p_\omega-\bar p_\omega\)$ with the slip law: 
\begin{equation}\label{eq_monitoring}
 \begin{aligned}
 &  q_\omega  \leq \bar q_\omega, \quad p_\omega \geq \bar p_\omega, \quad (q_\omega  -\bar q_\omega)(p_\omega - \bar p_\omega)=0, \\  
 & u^\g_\omega = F^\g(s^{\g}_{\s_\omega},s^{\g}_{\s_\omega},H_{\s_\omega}) = s^{\g}_{\s_\omega} C_0(s^{\g}_{\s_\omega}, H_{\s_\omega}) u^{\rm m}_\omega   + s^{\g}_{\s_\omega} U_{\rm d}(s^{\g}_{\s_\omega}, H_{\s_\omega}), \\
&  u^{\rm m}_\omega = u^\g_\omega + u^\l_\omega. 
\end{aligned}
\end{equation}

\section{Nonlinear solver} \label{sec:nonlinear}

Given the subset of present phases $Q_\s$, let us define the set of well unknowns at node $\s\in \mathcal{V}_\omega$
$$
X_\s = \(p_\s, \, T_\s, \, s_{\s}^\g, s_{\s}^\l, \, c^\alpha_\s, \alpha\in Q_\s, \, \tilde n_{\s,i}, i\in \wt{\mathcal{C}}_{Q_\s} \).
$$
The discrete nonlinear system to be solved at each time step is defined for each well by the set of unknowns 
$$
(X_\s)_{\s\in \mathcal{V}_\omega}, \, (Q_\s)_{\s\in \mathcal{V}_\omega},\, (u^\g_\welledge, u^\l_\welledge)_{\welledge \in \edges_\omega}, \, (u^\g_\omega, u^\l_\omega),
$$
and the set of equations \eqref{eq_well_discrete}-\eqref{eq_superficial_velocities}-\eqref{eq_monitoring} complemented by the flash equations of type $Q_\s = \mbox{\rm flash}(X _\s)$ for each node $\s\in \mathcal{V}_\omega$.
It is coupled fully implicitly to the reservoir system through the source terms $q^{r\rightarrow \omega}_{\s,i}$, $q^{r\rightarrow \omega}_{\s,e}$ for $\s \in \mathcal{V}_\omega$. Let us refer to Appendix \ref{sec:peaceman} for their expressions depending on the reservoir and well pressure, temperature, molar fractions and saturations at node $\s$. 

This system is solved using a Newton-Raphson algorithm based on an active set method for the subsets of present phases and for the well active constraints. This type of algorithm can be viewed as an active set formulation of a Newton-Min semi-smooth Newton algorithm (see e.g. \cite{Xing.ea:2017} and \cite{Beaude2019}).

The Jacobian system is assembled at each Newton iteration given the sets of present phases and the well active constraints which are updated after each Newton update. An important feature of the implementation is the reduction of the Jacobian system by elimination of the local thermodynamical closure laws \eqref{eq_thermo_discrete}, of the hydrodynamical equations \eqref{cons_momentum_discrete}-\eqref{eq_superficial_velocities} and of the well monitoring conditions \eqref{eq_monitoring}.
The elimination of the thermodynamical closure laws is also applied to the reservoir equations. 
The reduced linear system couples the component and energy conservation equations with a block structure of $\# \mathcal{C}+1$ equations and primary unknowns at each well and reservoir nodes. This elimination procedure detailed in the next subsections facilitates the assembly of the coupled linear system and reduce the cost of its resolution.

\subsection{Elimination of the thermodynamical closure laws}\label{sec:elimthermo}

Let us denote by $C_\s(X_\s)$ the thermodynamical closure equations \eqref{eq_thermo_discrete} at node $\s$.  
The elimination of the thermodynamical closure laws is based on a splitting of the unknowns $X_\s = (X_\s^p, X_\s^s)$
into $\# \mathcal{C}+1$ primary unknowns $X^p_\s$ and the remaining secondary unknowns $X^s_\s$. This splitting classically depends on the set of present phases $Q_\s$ and is such that ${\partial C_\$(X^p_\s,X^s_\s) \over \partial X_\s^s}$ is non singular. It results that, after Newton linearization, the secondary unknowns $dX^s_\s$ can be expressed as a linear function of $dX^p_\s$ and of the residual $C_\s(X_\s)$ (see \cite{Xing.ea:2017} for details in the case of the Coats' formulation of the reservoir model).

 \subsection{Elimination of the hydrodynamical equations}\label{sec:elimhydro}

For all well edge $\welledge = \s\s'$, given $X_\s, X_{\s'}, Q_\s, Q_{\s'}$, the hydrodynamical module computes the phase superficial velocities
$$
u_\welledge^\alpha  \mbox{ for } \alpha \in Q_\s \cup Q_{\s'},
$$
and their Newton linearization accounting for the elimination of the thermodynamical closure laws, such that
\begin{equation}\label{elim:vs}
d u_\welledge^\alpha = A^\alpha_{\welledge,\s} dX^p_\s + A^\alpha_{\welledge,\s'} dX^p_{\s'} + B^\alpha_\welledge \mbox{ for } \alpha \in Q_\s \cup Q_{\s'}.  
\end{equation}
This computation is detailed in the next two paragraphs. 
Note that, from the properties of the numerical flux function $F_\welledge^\g$,  we have $u_\welledge^\alpha =0$ for  $\alpha \not\in Q_\s \cup Q_{\s'}$, hence this case does not need to be considered. 

\subsubsection{Mixture velocity}

Let us define the difference of potential at the well edge $\welledge=\s\s'$ by
$$
\Delta_\welledge \Phi =  p_{{\s}} - p_{{\s'}} + \bar \rho^{\rm m}_\welledge g (z_{{\s}} - z_{{\s'}}).  
$$
Considering the wall friction law  \eqref{eq_darcy_forchheimer}, the momentum equation \eqref{cons_momentum_discrete} can be solved for the mixture velocity as a function of $\Delta_\welledge \Phi$ depending on the sign of 
$\Delta_\welledge \Phi$ as follows 

\begin{equation}
  \label{eq_um}
u_\welledge^{\rm m} =  - sign(\Delta_\welledge \Phi) {(\alpha_c - \alpha_b)  \over 2 \alpha_a},
\end{equation}
with
$$
\alpha_a = |\welledge| f_q  { \bar \rho^{\rm m} \over 4 r_\welledge },\quad \alpha_b = 8 |\welledge|{ \bar \mu^{\rm m} \over (r_\welledge)^2 }, \quad \alpha_c = \sqrt{ (\alpha_b)^2 + 4  |\Delta_\welledge \Phi| \alpha_a },  
$$
and the definition \eqref{eq_mum_rhom} of the mean specific density $\bar \rho^m_\welledge$ and viscosity $\bar \mu^m_\welledge$. 

Note that for more general wall friction laws, the elimination of the mixture velocity is not always possible at the non linear level. In that case it is done at the linear level after Newton linearization taking into account the momentum equation residual. 
The mixture velocity $u^{\rm m}_\welledge$ is computed from \eqref{eq_um} as a function of the entries
$$
\Delta_\welledge \Phi, \,\, \bar \rho^{\rm m}_\welledge, \,\, \bar\mu^{\rm m}_\welledge, 
$$
as well as its derivative w.r.t. the 3 entries.

Then, using the two-point monotone flux function $F^g_\welledge$, we can compute both superficial velocities as detailed below. 

\subsubsection{Superficial velocities}
For $Q_\s \cup Q_{\s'} = \{\alpha\}\in \mathcal{P}$, $\{\beta\} = \mathcal{P}\setminus \{\alpha\}$, we simply have
$$
u^\alpha_\welledge = u^{\rm m}_\welledge,\qquad u^\beta_\welledge = 0, 
$$
with $u^{\rm m}_\welledge$ defined by \eqref{eq_um}. \\

For $Q_\s \cup Q_{\s'} = \mathcal{P}$, we first compute $u^{\rm m}_\welledge$ from \eqref{eq_um}, and
then, the gas superficial velocity is obtained from the two-point monotone flux $F^\g_\welledge$.
Using the model detailed in Appendix \ref{sec_monotoneflux}, the entries of  $F^\g_\welledge$ are 
$$
s^\g_{\s'},\,\, s^\g_{\s},\,\, \bar \rho_\welledge^\l, \,\, \bar \rho_\welledge^\g,\,\, \bar \sigma_{\g\l,\welledge}, \,\, u^{\rm m}_\welledge. 
$$
The flux function $F^\g_\welledge$ computes $u^\g_\welledge$ as well as its derivatives w.r.t. the 6 entries.
Then the liquid superficial velocity is given by
$$
u^\l_\welledge = u^{\rm m}_\welledge - u^\g_\welledge. 
$$
Combining the derivatives of the mixture velocities and superficial velocities w.r.t. to the above entries with the derivatives of the thermodynamical variables, we can compute \eqref{elim:vs}. 

\subsection{Elimination of the monitoring conditions}\label{sec:elimmonitoring}

The complementary constraints are solved using an active set method meaning that we impose either a fixed pressure $p_\omega = \bar p_\omega$ with $q_\omega  \leq \bar q_\omega$ or a fixed molar flow rate $q_\omega  = \bar q_\omega$ with $p_\omega \geq \bar p_\omega$. The active constraint is updated at each Newton iteration according to the remaining inequality constraint. 

The elimination of the head node superficial velocities is performed as follows. 
In case of active flow rate constraint, the well total molar flow rate is prescribed with 
\begin{equation}
  \label{eq_qfixed}
 q_\omega =  \bar q_\omega, 
\end{equation}
 In case of active pressure constraint, it is obtained by the sum over $i$ of the molar conservation equations at node $\s_\omega$
 \begin{equation}
    \label{eq_pfixed}
q_\omega =  \dsp - \sum_{i\in\mathcal{C}} {n_{\s_\omega,i}- n_{\s_\omega,i}^{n-1}\over \Delta t^n } + \dsp \sum_{i\in\mathcal{C}} q^{r\rightarrow \omega}_{\s_\omega,i} + \sum_{\welledge\in \edges_{\s_\omega}^\omega}  \sum_{\alpha\in \mathcal{P}} \kappa_{\welledge,\s_\omega} |S^\omega_\welledge|\zeta^\alpha_\welledge u^\alpha_\welledge. 
 \end{equation}
Note that this equation must hence be substituted by the equation $p_\omega = \bar p_\omega$ in the system of primary equations and unknowns. 
 Using the total molar flow rate $q_\omega$ computed from \eqref{eq_qfixed} or \eqref{eq_pfixed}, we can eliminate the head node mixture velocity $u^{\rm m}_\omega$ as follows. Using \eqref{eq_wellpwqw} and \eqref{eq_monitoring}, we set
 $$
    u^{\rm m}_\omega = { q_\omega \over |S^\omega_{\s_\omega}| \zeta^\alpha_{\s_\omega} } \mbox{ if } Q_{\s_\omega}  = \{\alpha\},  
$$
 and,  if $Q_{\s_\omega} = \mathcal{P}$, setting  $C_{0,\s_\omega} = C_0(s^{\g}_{\s_\omega}, H_{\s_\omega})$, $U_{{\rm d},\s_\omega} = U_{\rm d}(s^{\g}_{\s_\omega}, H_{\s_\omega})$, we obtain the equation
 $$
 {q_\omega \over |S^\omega_{\s_\omega}|} = \zeta^\g_{\s_\omega} \(s^{\g}_{\s_\omega} C_{0,{\s_\omega}} u^{\rm m}_\omega + s^{\g}_{\s_\omega} U_{{\rm d},{\s_\omega}} \)
 + \zeta^\l_{\s_\omega}  \( (1- s^{\g}_{\s_\omega} C_{0,{\s_\omega}}) u^{\rm m}_\omega - s^{\g}_{\s_\omega} U_{{\rm d},{\s_\omega}}\).  
 $$
providing 
$$
 u^{\rm m}_\omega ={  {q_\omega \over |S^\omega_{\s_\omega}|} +(\zeta^\l_{\s_\omega} - \zeta^\g_{\s_\omega}) s^{\g}_{\s_\omega} U_{{\rm d},{\s_\omega}}
   \over \zeta^\g_{\s_\omega}  s^{\g}_{\s_\omega} C_{0,{\s_\omega}} + \zeta^\l_{\s_\omega}  (1-s^{\g}_{\s_\omega} C_{0,{\s_\omega}}) }. 
$$

 Once $u^{\rm m}_\omega$ is known we can compute the head node superficial velocities $u^\alpha_\omega$ from \eqref{eq_monitoring}.

 \section{Numerical experiments} \label{sec:testcases}

 The objectives of this numerical Section is to validate the numerical model and investigate its ability to simulate both well opening transient test cases without and with cross flow and  fully coupled test cases at the reservoir time scale. The well stand alone test cases are presented in Subsection  \ref{sec:standalone} and a fully coupled test case is considered in Subsection \ref{sec:fullycoupled}  where it is compared for validation to a simpler single implicit unknown well model.

%The advantages of the modeling approach described in the previous section are illustrated by four test cases. The three first cases concern stand alone well model. The first numerical experiment is a two-phase incompressible immiscible model, which reduces to a scalar hyperbolic Buckley Leverett equation. The nonlinear solution obtained is compared with that of another model using a different discretization. The two next examples are dedicated to the study of the production of a multi-branch well without or with cross-flow. The last test case is devoted to coupling the well with the geothermal reservoir. Considering the production of a vertical well, we solve the coupled problem of two-phase flow in the porous medium and two-phase flow in the well, and compare the result obtained with a simpler production well model with a single unknown per well.

\subsection{Stand alone well model}\label{sec:standalone}

 \subsubsection{Validation on a Buckley Leverett solution}\label{sec:buckleyleverett}

 Let us consider a two component $\mathcal{C} = \{1,2\}$, two-phase $\mathcal{P} = \{\l,\g\}$, immiscible, incompressible and isothermal model obtained by setting in the general compositional framework $\mathcal{C}^\l = \{1\}$ and  $\mathcal{C}^\g = \{2\}$. 
The well is vertical of diameter $0.1$ m and centerline $(0,H)$ with $z=0$ corresponding to the leaf node and $z=H=100$ m to the head node.
 The molar and specific densities as well as the dynamic viscosities are considered constant.

\begin{figure}[H]
\begin{center}
\includegraphics[width=0.35\textwidth]{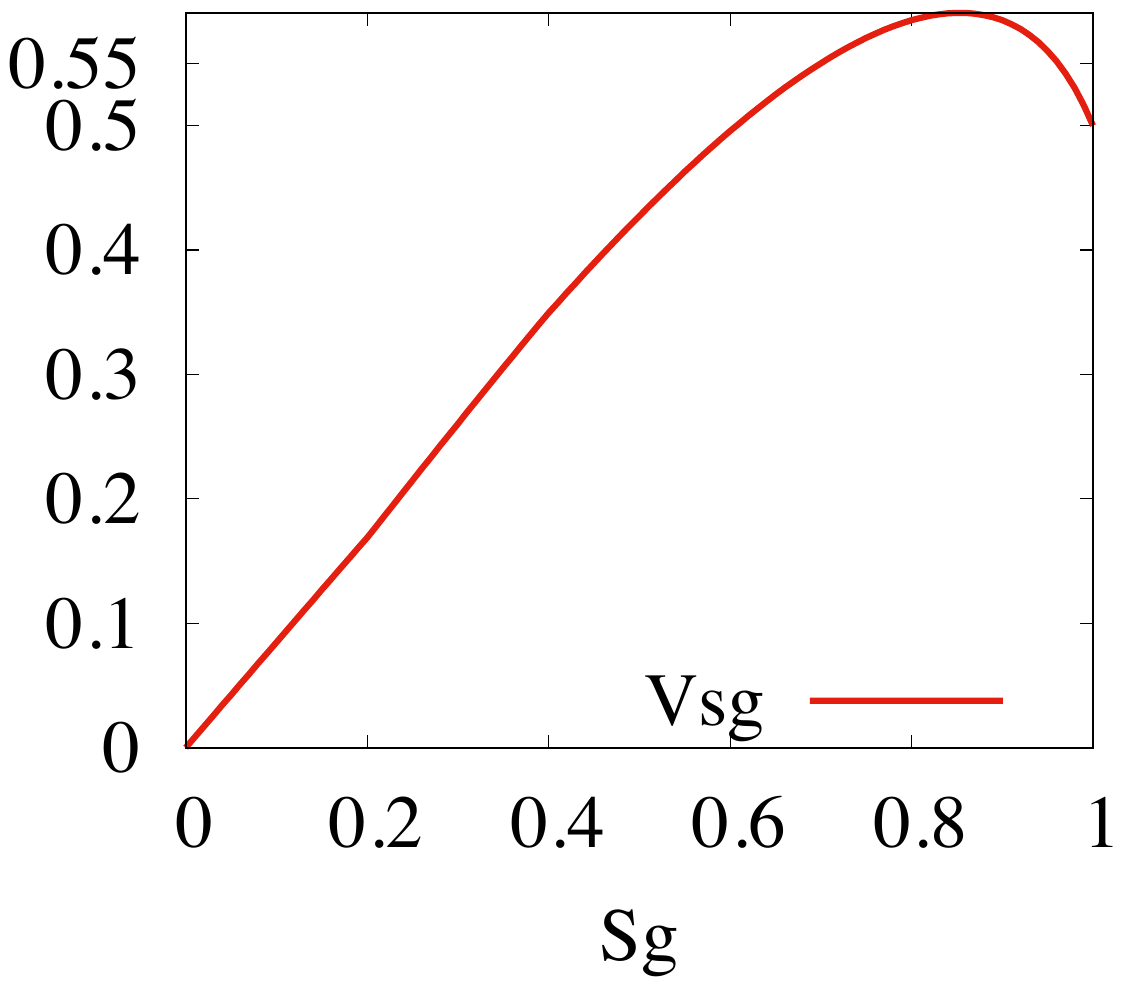}
\includegraphics[width=0.35\textwidth]{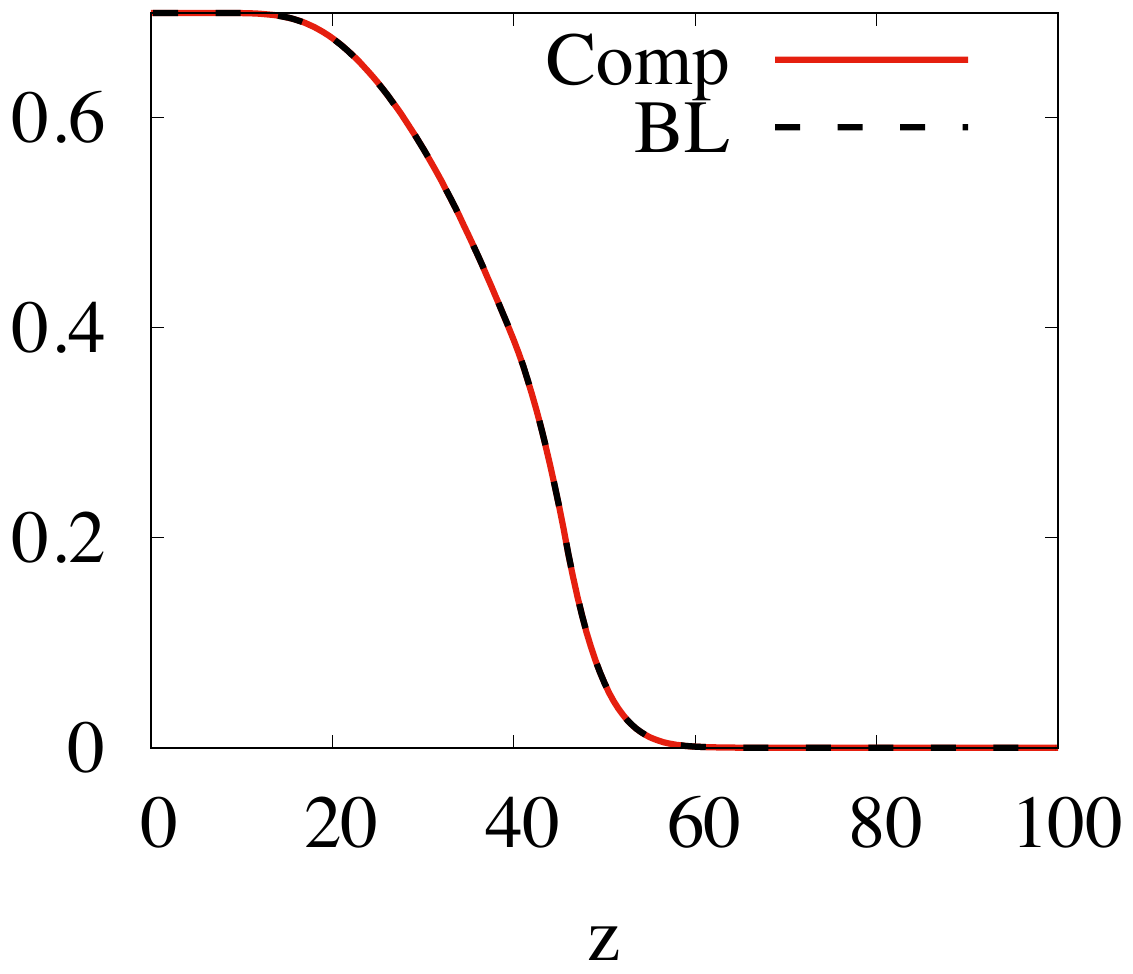}     
\caption{(Left): gas superficial velocity as a function of the gas saturation (flux function) at fixed mixture velocity $u^{\rm m} = 0.5$ m.s$^{-1}$ for the Buckley Leverett test case. (Right): Comparison between the Buckley Leverett (BL) and compositional (Comp) gas saturations at final time on the Buckley Leverett test case}
\label{fig_test1_Vsg_twocodes}
\end{center}
\end{figure}

\begin{figure}[H]
\begin{center}
  \includegraphics[width=0.35\textwidth]{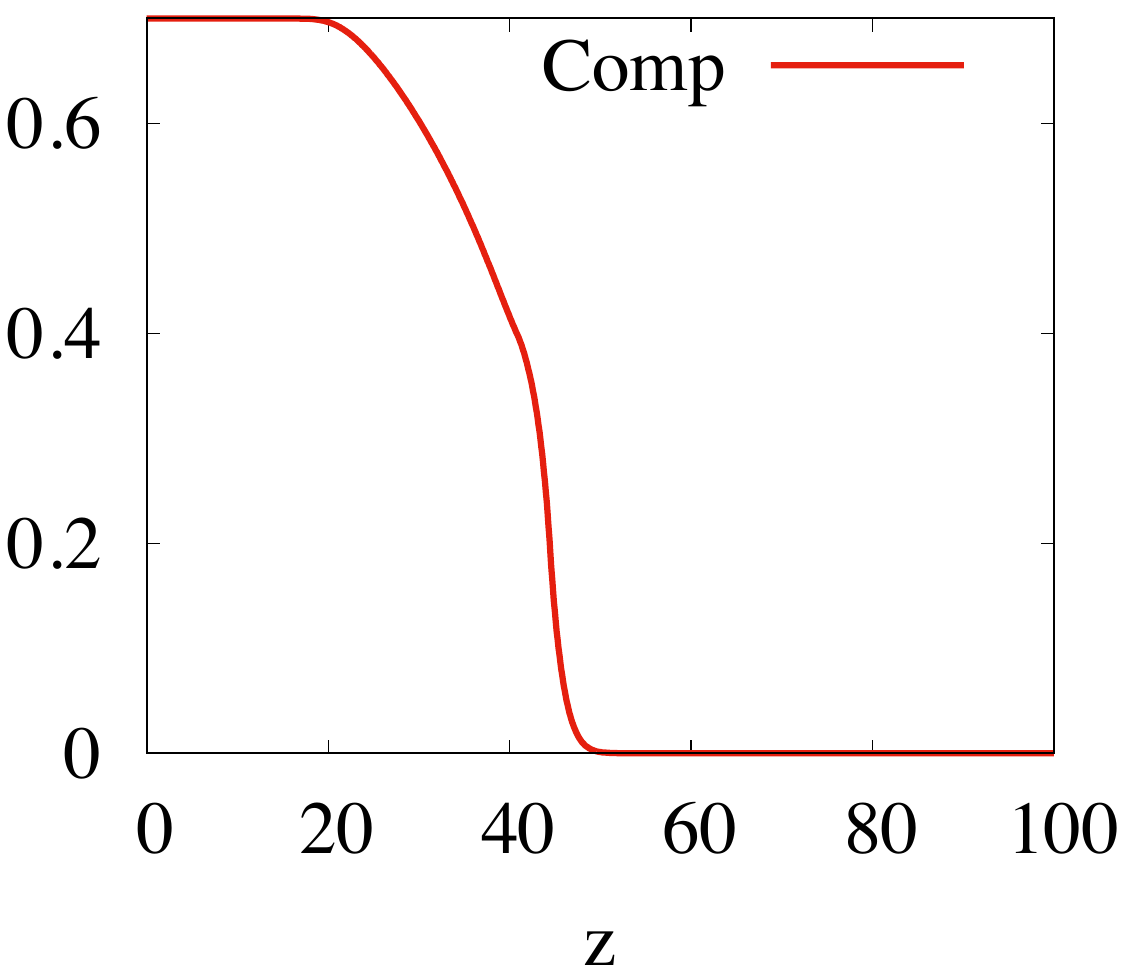}
\includegraphics[width=0.35\textwidth]{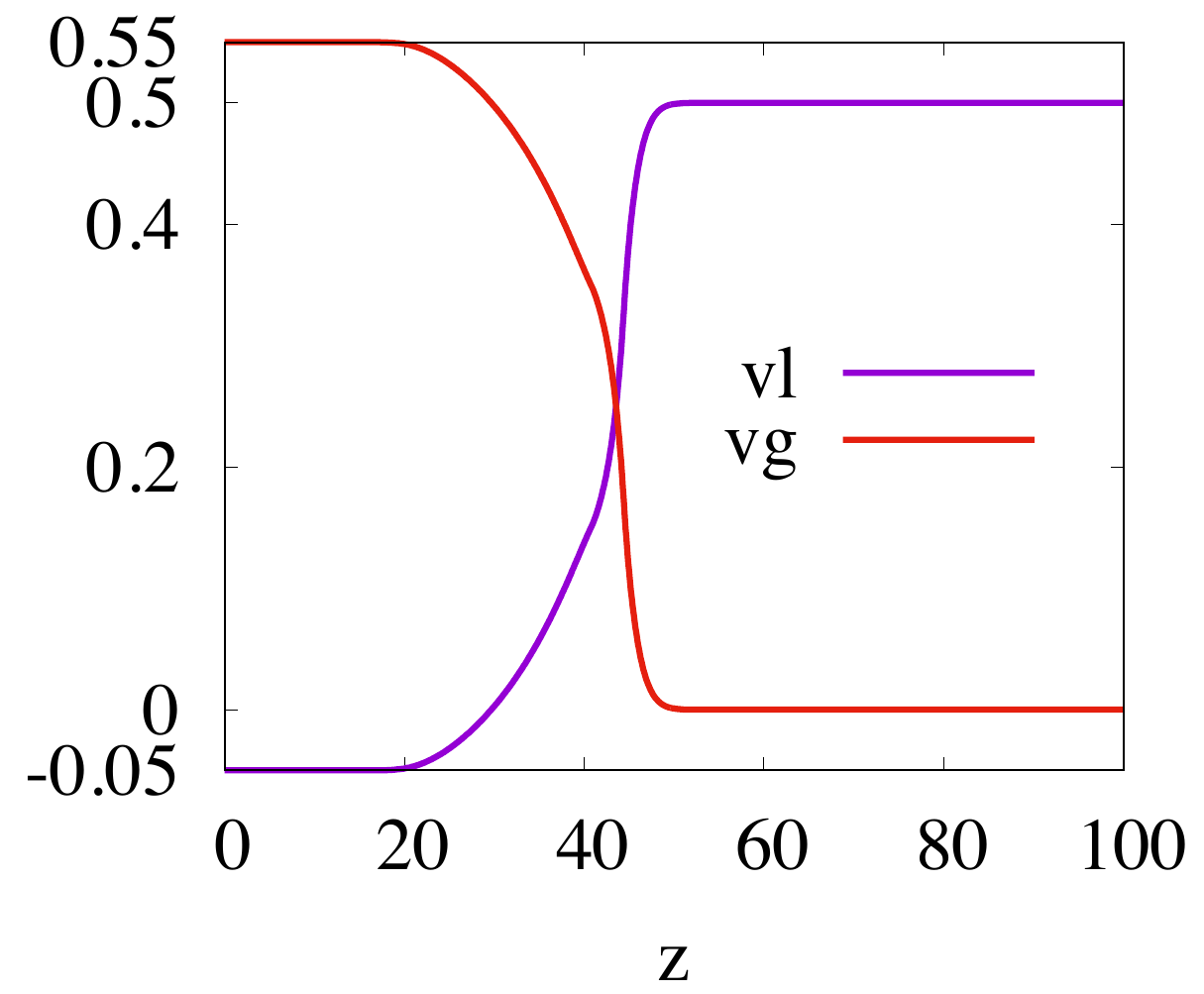}     
\caption{ (Left): Compositional gas saturation at final time using the fine mesh with $1000$ cells and $1000$ time steps on the Buckley Leverett test case. (Right): liquid and gas superficial velocities obtained with the compositional model at final time using the fine mesh with $1000$ cells and $1000$ time steps on the Buckley Leverett test case.  }
\label{fig_test1_sg_vlg}
\end{center}
\end{figure}

 Let us set
 $$
 u^\g(s^\g) = s^\g U_{\rm d}(s^\g)  +   s^\g C_0(s^\g) u^{\rm m},
 $$
 and
 $$
 u^\l(s^\g) = u^{\rm m} - u^\g(s^\g). 
 $$
 The volume conservation equations for each phase write
$$
\left\{\begin{array}{r@{\,\,}c@{\,\,}l}
& \partial_t s^\g + \partial_z u^\g(s^\g) = 0,\\
& \partial_t s^\l + \partial_z u^\l(s^\g) = 0.   
\end{array}\right.
$$ 
Summing the above equations we obtain $\partial_z u_m = 0$ meaning that the mixture velocity depends only on time.
It will be fixed by the input boundary condition at the leaf node fixing the mixture velocity to $u^{\rm m} = 0.5$ m.s$^{-1}$. 
Thus, we obtain the scalar hyperbolic Buckely Leverett equation
$$
\partial_t s^\g + \partial_z u^\g(s^\g) = 0,
$$
with the flux function $u^\g(s^\g)$ at fixed mixture velocity.
In this test case we compare the solution of the model developped in this work to the one obtained by solving numerically the scalar hyperbolic equation. The same Euler implicit time discretization is used for both implementations. The scalar hyperbolic discretization is based on a cell centered Finite Volume (FV) scheme while the compositional model uses the node centered FV scheme described in this work. 
The data set is defined by $\rho^g = 4$ Kg.m$^{-3}$, $\rho^l = 1000$ Kg.m$^{-3}$, by the input gas superficial velocity at the leaf node $u^\g = 0.55$ m.s$^{-1}$ and the output liquid velocity at the leaf node $u^\l = -0.05$ m.s$^{-1}$. 
The DFM model \cite{shi2005} (see Appendix \ref{sec_monotoneflux}) is used with the parameters $A = 1.2$, $B=0.3$, $a_1 = 0.2$, $a_2=0.4$, $K_u = 1.5$ and $\sigma_{\g\l}$ is fixed to $71.97 ~10^{-3}$ (see the resulting flux function on the left figure \ref{fig_test1_Vsg_twocodes}). 
The solutions are computed with both models using a uniform mesh of $200$ cells, a uniform time stepping with $200$ time steps and  the simulation time $t_F=50$ s.
The right figure \ref{fig_test1_Vsg_twocodes} shows that both models provide basically the same gas saturation plot at final time $t_F$. 
We also plot in Figure \ref{fig_test1_sg_vlg} the fine mesh gas saturation (Left) and superficial velocities (Right) at final time obtained with the compositional model using $1000$ cells and $1000$ time steps. 
The gas moves up at a higher velocity than the mixture velocity as a result of the drift velocity allowing the liquid to go out of the well at the leaf node. 

\subsubsection{Thermal test case with a chair shaped production well}\label{sec:chairwell.test}

\begin{figure}[H]
\begin{center}
\includegraphics[width=0.1\textwidth]{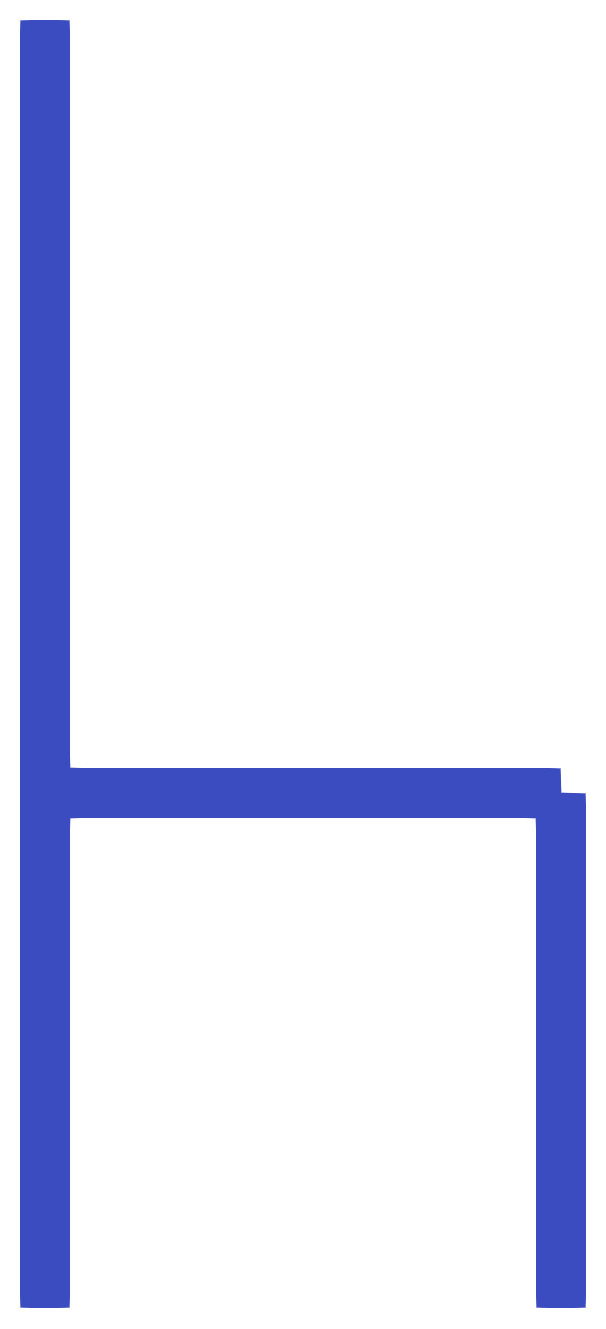}   
\caption{Chair shaped well with one junction and four branches of sizes 600 m (upper vertical branch), 400 m (lower left vertical branch), 400 m (lower vertical branch) and 400 m (horizontal branch).}
\label{fig_test2_chairwell}
\end{center}
\end{figure}

We consider a single H$_2$O component liquid vapor thermal flow along the chair shaped production well illustrated in Figure \ref{fig_test2_chairwell}.
The internal energy, mass density and viscosity of H$_2$O  in the liquid and gas phases  are defined by analytical laws as functions of the pressure and temperature (refer \cite[Section 4.1]{ACJLM23} for details). The vapour pressure $P_{sat}(T)$ is given in Pa by 
$$
p_{\rm sat}(T) = 10^{-3} (T-273)^4.  
$$
The DFM model \cite{shi2005} is used with the parameters $A = 1.2$, $B=0.3$, $a_1 = 0.2$, $a_2=0.4$, $K_u = 1.5$ and $\sigma_{\g\l}$ is fixed to $71.97 ~10^{-3}$ (see Appendix \ref{sec_monotoneflux}). The Darcy-Forccheimer friction law \eqref{eq_darcy_forchheimer} is set up with the friction parameter $f_q = 6~10^{-2}$. The well radius is fixed to $0.05$ m, the thermal conductivity is considered constant for simplicity with $\lambda = 2$.

The well is monitored with the minimum well head pressure set to $\bar p_\omega = 5 ~10^5$ Pa and the maximum well mass flow rate set to $\bar q_\omega = 15$ Kg.s$^{-1}$. The initial temperature is set to $350$ K along the well and the initial pressure is hydrostatic at liquid state with pressure set to $5 ~10^5$ Pa at the head node.

Feed zones are modeled at each of the two leaf nodes at high temperature and high pressure in liquid state. At both leaves the reservoir pressure is fixed to $1.1~10^7$ Pa, the temperature to $520$ K, the 
Darcy well index to $WI^D = 10^{-12}$ m and the Fourier well index to $WI^F = 100$ J.s$^{-1}$.K$^{-1}$. 

The well is meshed using a uniform discretization of each of the four branches with $40$ edges for both bottom vertical branches, $40$ edges for the horizontal branch and $60$ edges for the top vertical branch. The simulation time is fixed to $t_F=2000$ s at which the stationnary state is basically reached, and the time stepping is set up with an initial time step of $10$ s and a maximum time step of $40$ s.
The stopping criteria of the Newton nonlinear solver is fixed to either $10^{-8}$ on the relative residual $l^1$ norm or to $10^{-10}$ on the $l^\infty$  norm of the Newton increment $ds^g + {dp \over 10^5} + {dT \over 100}$.
The time step is multiplied by the factor $1.1$ until it reaches the maximum time step in case of Newton convergence in less than $50$ iterations and restarted with a twice smaller time step otherwise. 
Using this setting, the simulation runs in $58$ time steps with no time step failure and a total number of $438$ Newton iterations.  

\begin{figure}[H]
\begin{center}
  \includegraphics[width=0.2\textwidth]{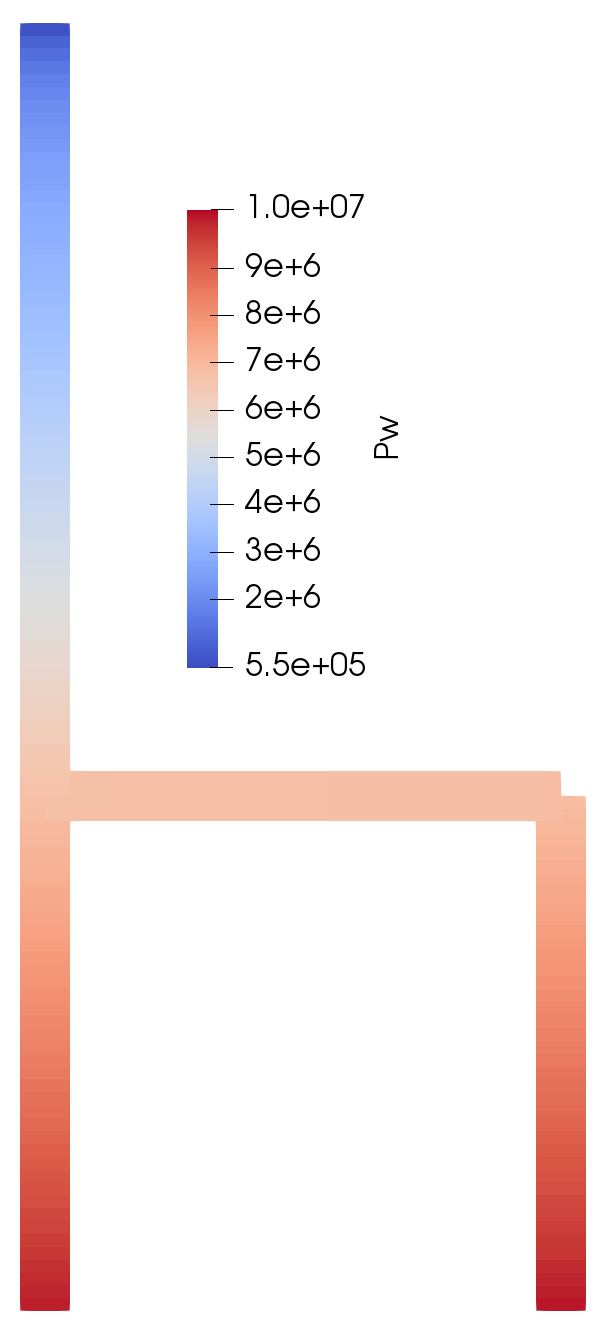}
  \includegraphics[width=0.2\textwidth]{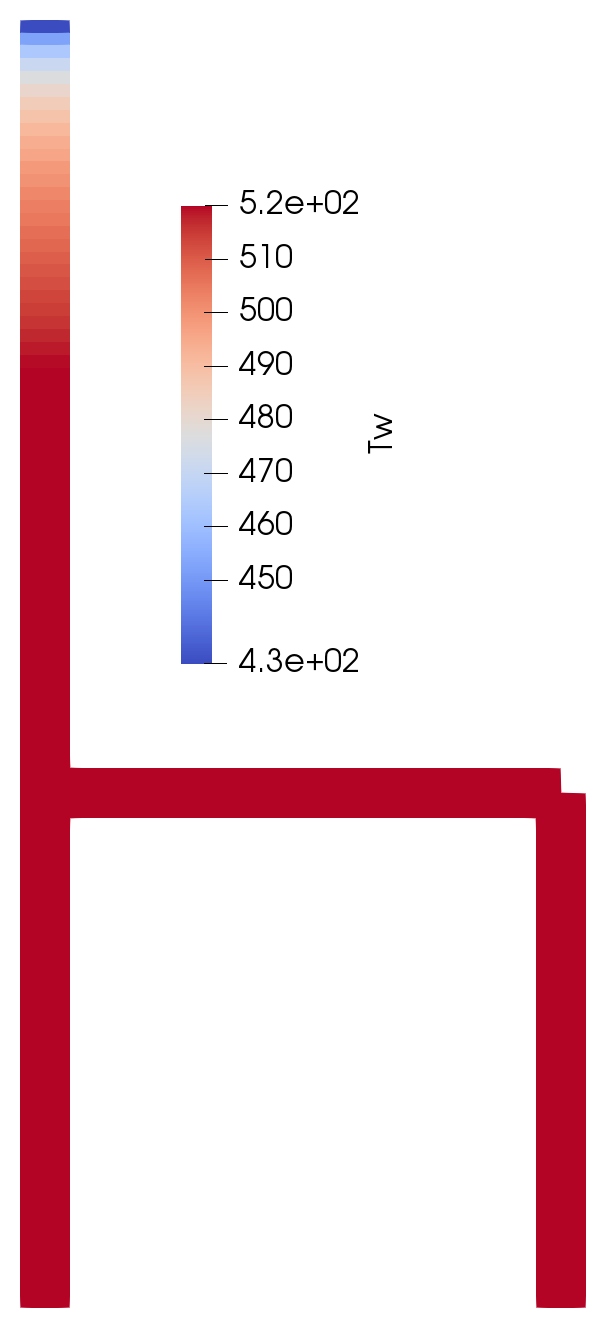}
  \includegraphics[width=0.2\textwidth]{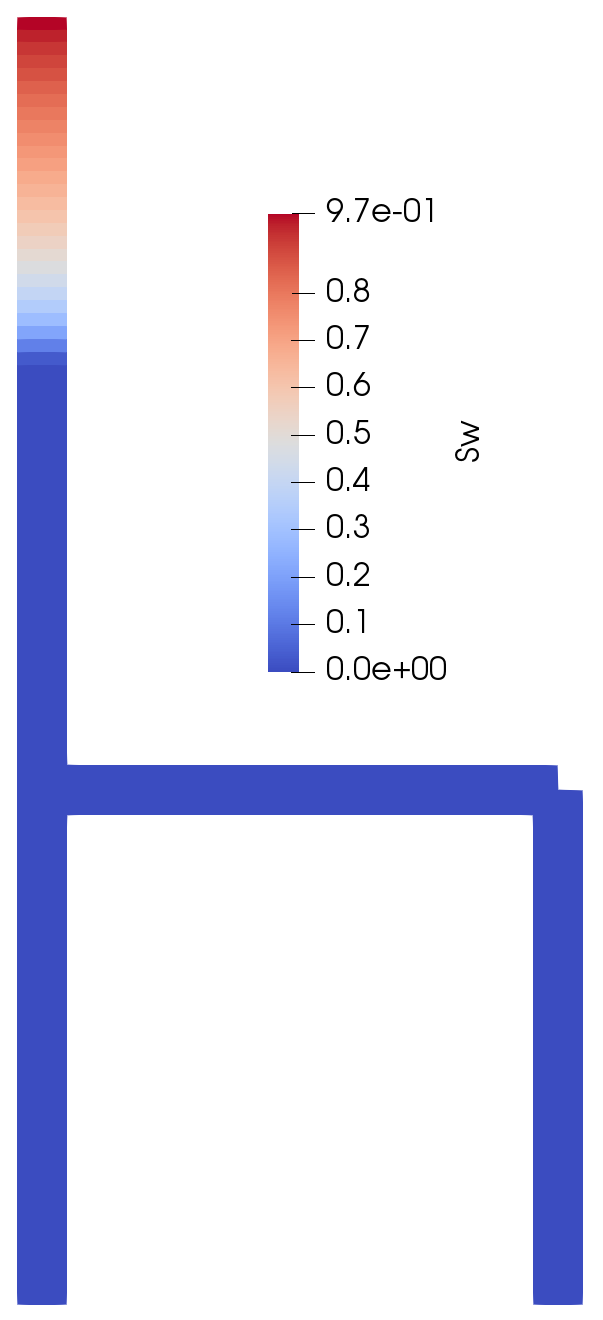}   
\caption{Pressure (Pa), temperature (K), gas saturation solutions at final time for the chair shaped thermal well test case.}
\label{fig_test2_PTS_paraview}
\end{center}
\end{figure}

\begin{figure}[H]
\begin{center}
  \includegraphics[width=0.2\textwidth]{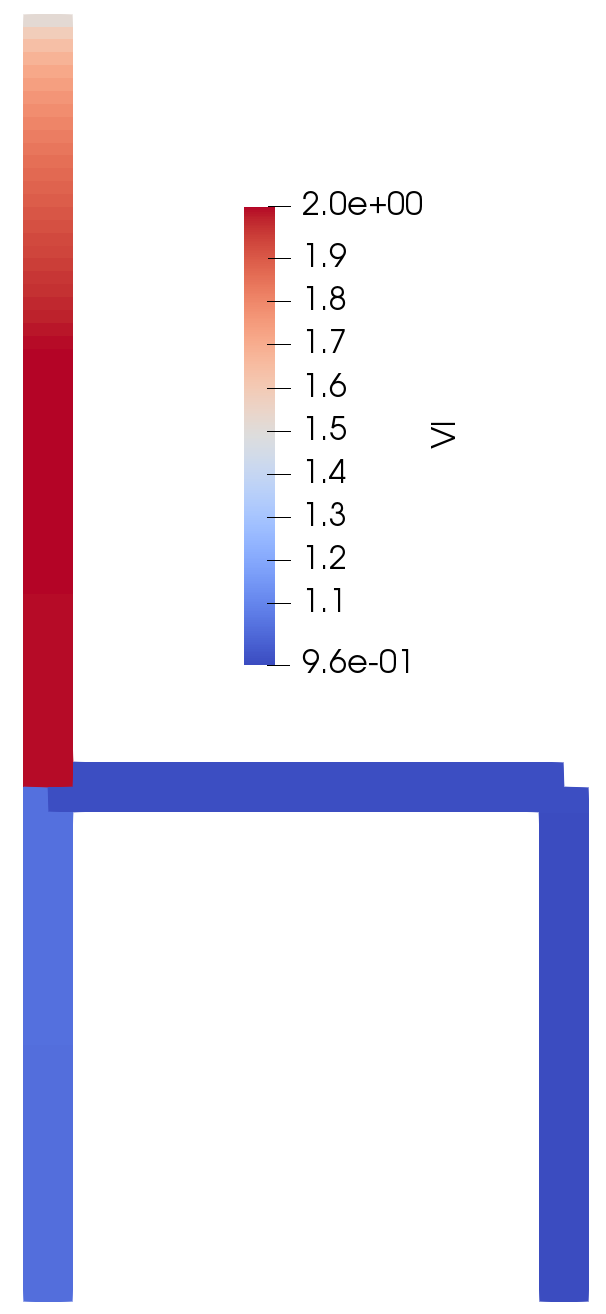}
  \includegraphics[width=0.2\textwidth]{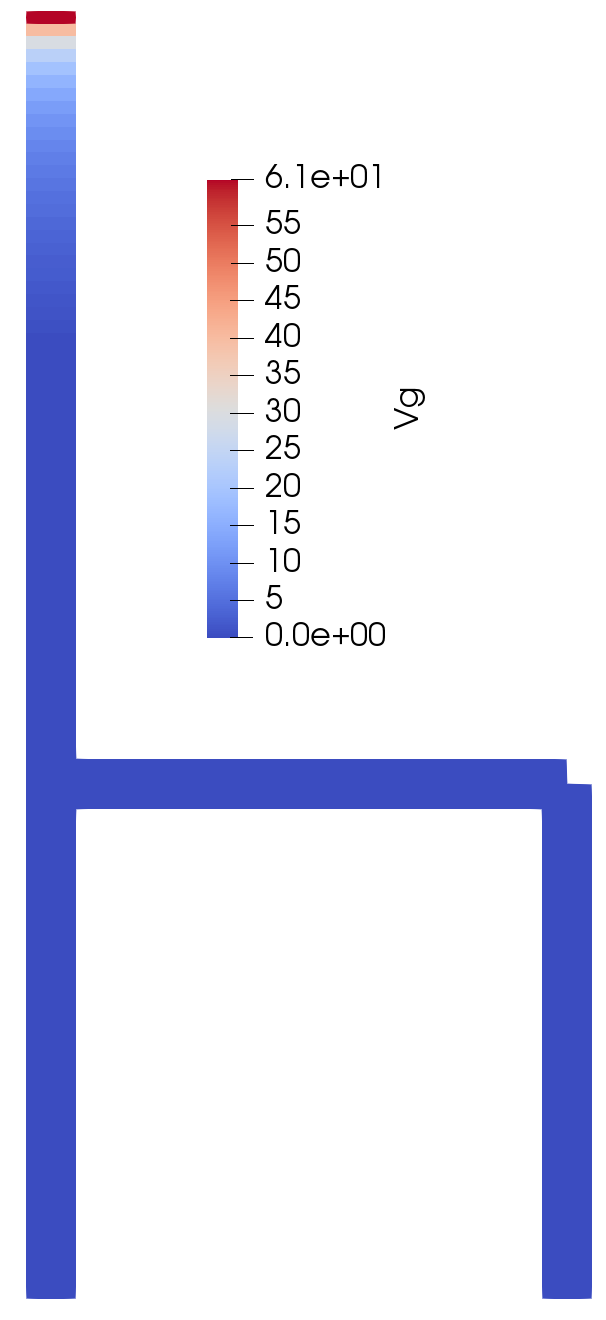}
\caption{Liquid and gas superficial velocities (m.s$^{-1}$) at final time for the chair shaped thermal well test case.}
\label{fig_test2_Vlg_paraview}
\end{center}
\end{figure}

\begin{figure}[H]
\begin{center}
  \includegraphics[width=0.4\textwidth]{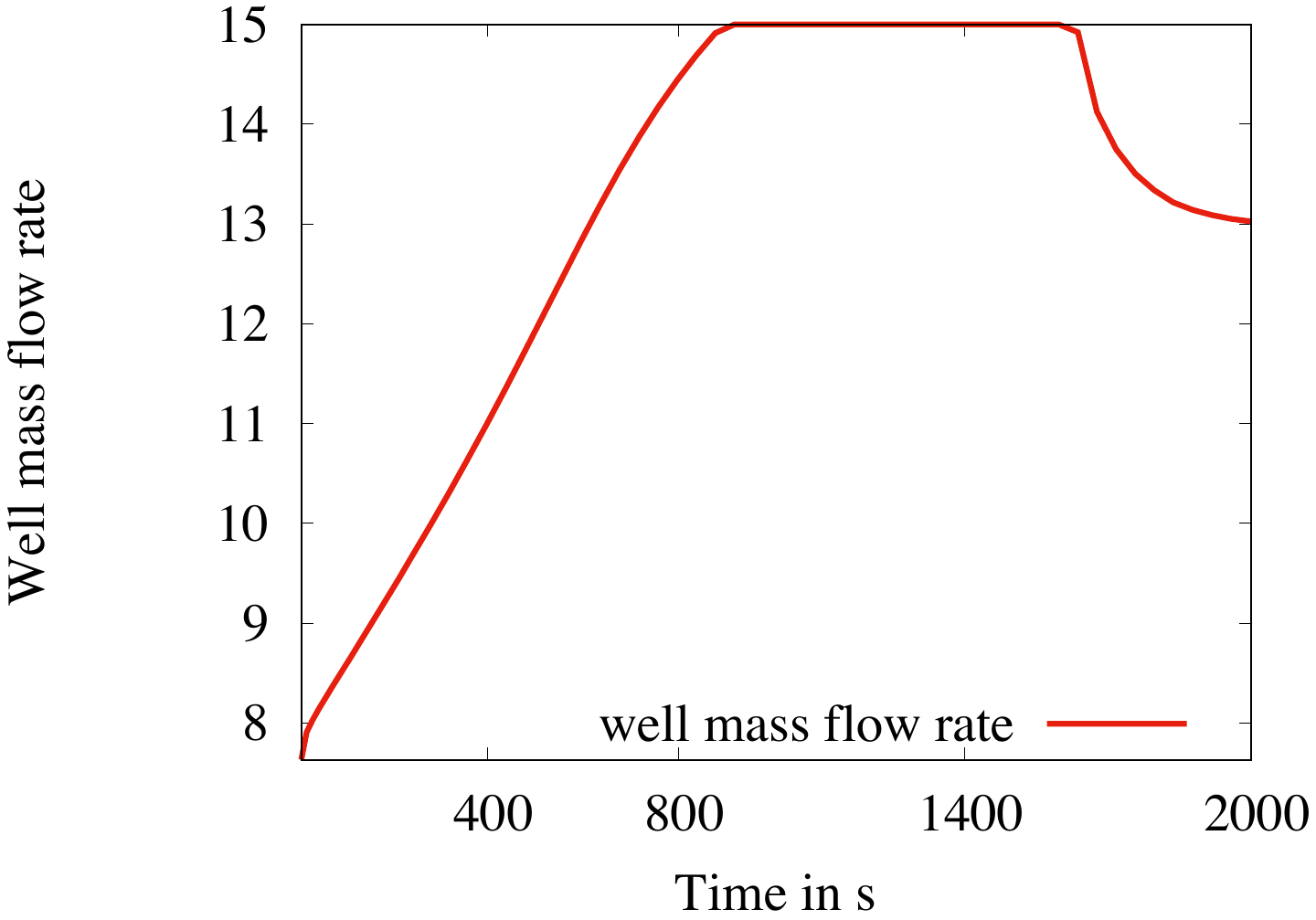}
  \includegraphics[width=0.4\textwidth]{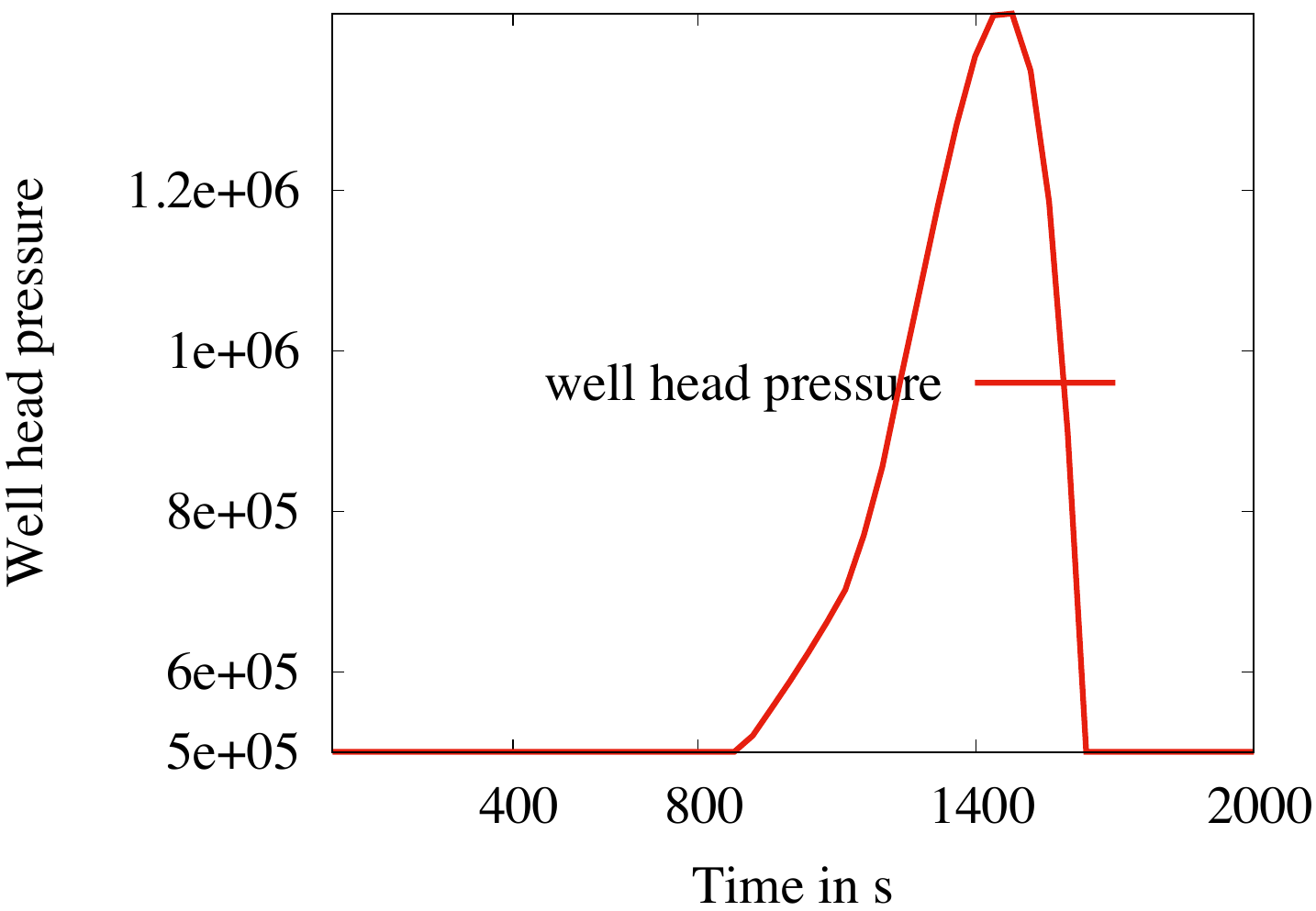}
  \includegraphics[width=0.4\textwidth]{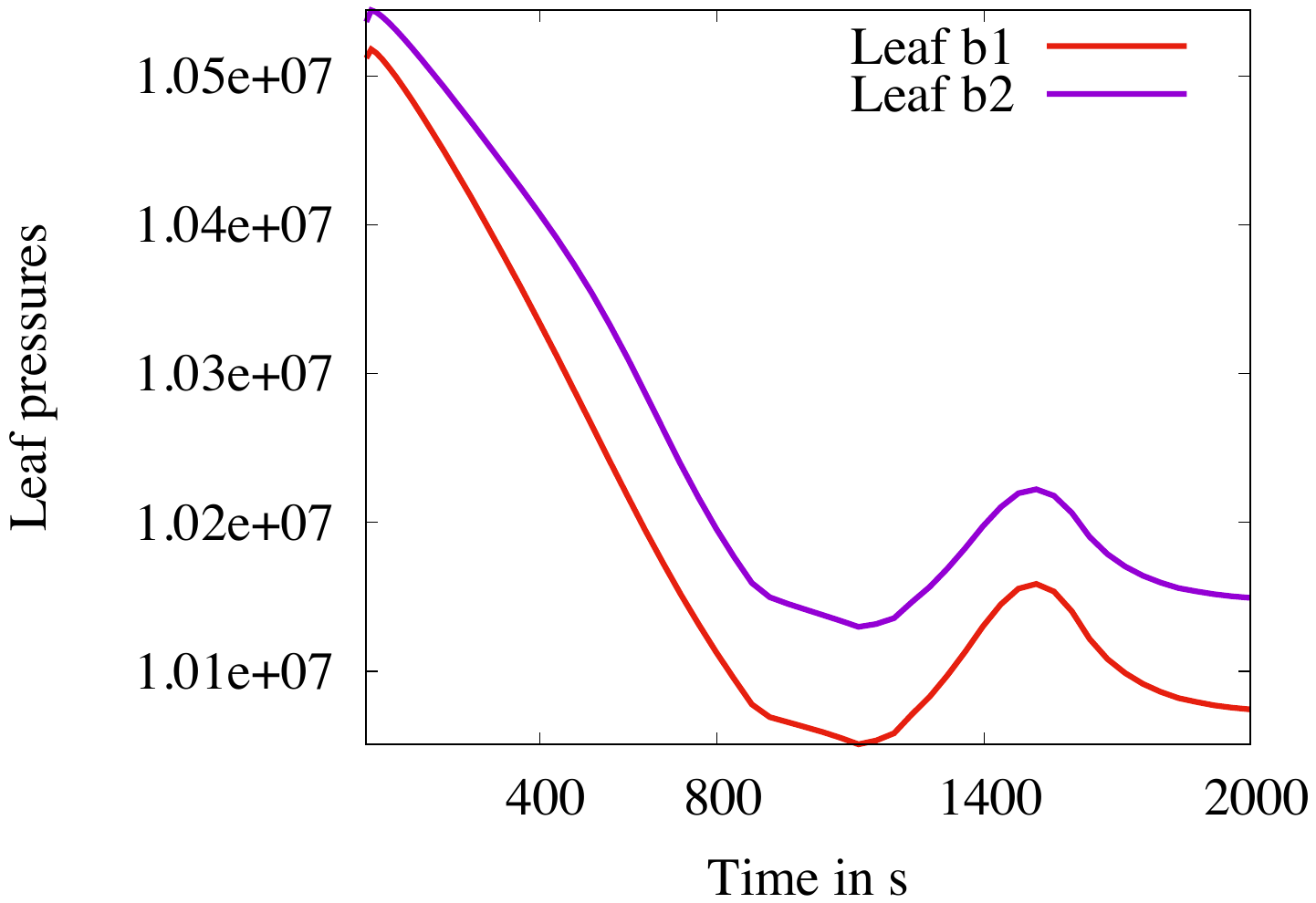}
   \includegraphics[width=0.4\textwidth]{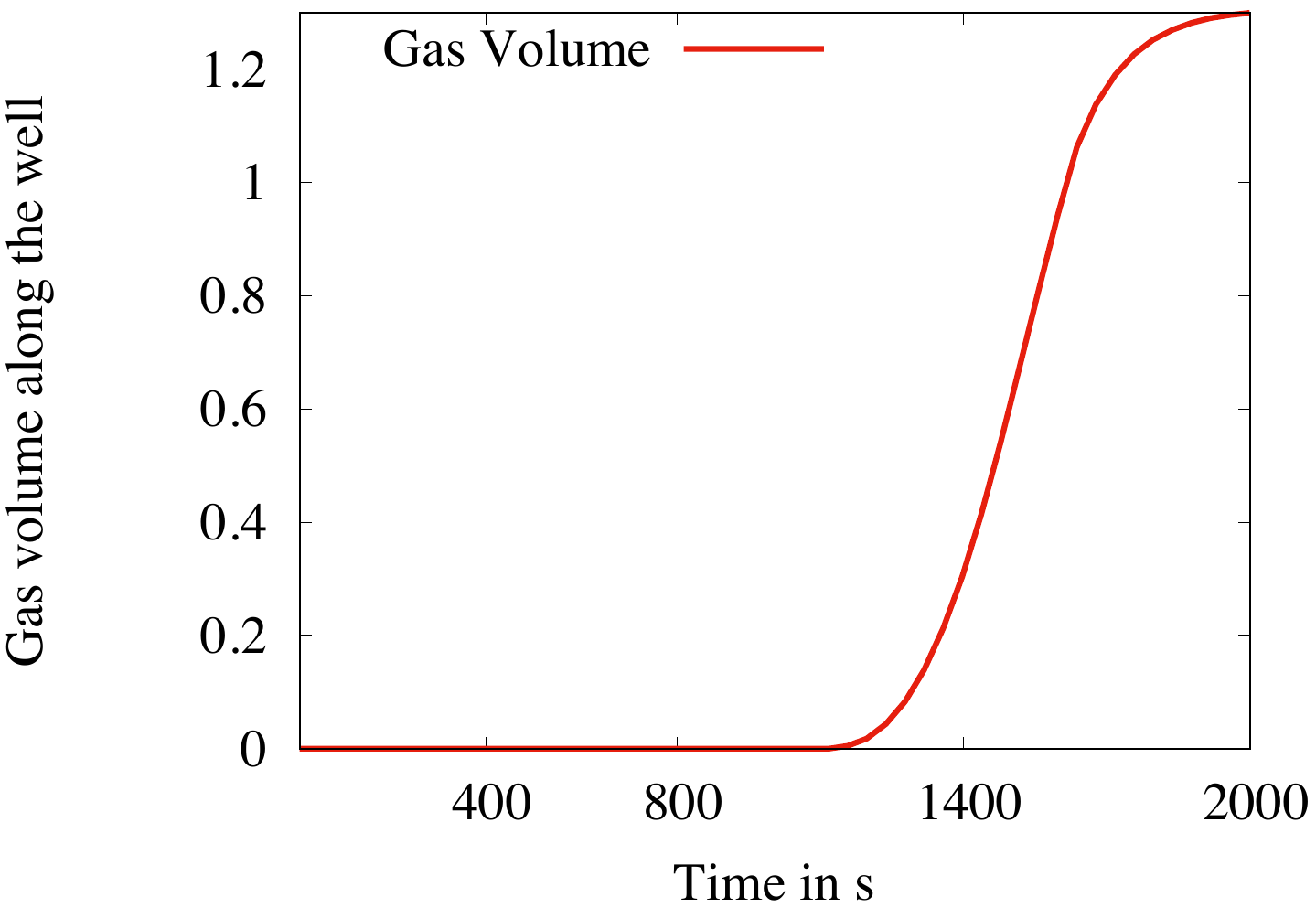}   
\caption{Well mass flow rate, well head pressure, leaf pressures, and gas volume in the well as a function of time for the chair shaped thermal well test case.}
\label{fig_test2_histories}
\end{center}
\end{figure}

Figures \ref{fig_test2_PTS_paraview} and \ref{fig_test2_Vlg_paraview} exhibit the final pressure, temperature, gas saturation and superficial velocities. The rise of the hot temperature front induces the appearance of the gas phase at the top of the well starting at roughly $t=1200$ s (see Figure \ref{fig_test2_histories}). Figure \ref{fig_test2_histories} plots the time histories of the well pressure $\bar p_\omega$, the mass flow rate $\bar q_\omega$, the leaf pressures and the gas volume inside the well.  The well starts to produce at the fixed minimum well head pressure of $5 ~10^5$ Pa with a rising mass flow rate. This rising flow rate is induced by the temperature increase along the well which reduces the weight of the liquid column. This can be checked in the leaf pressures plot showing the decrease of the leaf pressures in the first part of the simulation at fixed minimum well pressure. Then, the well flow rate reaches its maximum value of 15 Kg.$s^{-1}$ and the well operates at fixed maximun mass flow rate with a rising well head pressure which speeds up when the gas phase appears. At around $t=1400$ s, the well head pressure starts to decrease rapidly as a result of the pressure drop induces by the high gas velocity, until the well operates again at fixed minimum well head pressure of $5 ~10^5$ Pa. We can also notice, in the leaf pressures plot, the gap between the two leaf pressures which results from the additional wall friction along the horizontal branch.

%Figures \ref{fig_test2_PTS} and \ref{fig_test2_Vlg} exhibit the plots of the pressure, temperature, gas saturation and superficial velocities at different times, showing separately the vertical left branches (b1) in full lines and the right vertical + horizontal branches (b2) in dots lines. It provides a more accurate visualization than the paraview figures and clearly exhibits the discontinuities of $T$ (due to the convection dominated regime) and of $u^\l$ at the junction.  
%\begin{figure}[H]
%\begin{center}
%  \includegraphics[width=0.4\textwidth]{TestThermal/Pw.pdf}
%  \includegraphics[width=0.4\textwidth]{TestThermal/Tw.pdf}
%  \includegraphics[width=0.4\textwidth]{TestThermal/Sgw.pdf}   
%\caption{Pressure, temperature and gas saturation at different times along the left vertical left branches (b1) and the right vertical + horizontal branches (b2) for the chair shaped thermal well test case.}
%\label{fig_test2_PTS}
%\end{center}
%\end{figure}

%\begin{figure}[H]
%\begin{center}
%  \includegraphics[width=0.4\textwidth]{TestThermal/Vlw.pdf}
%  \includegraphics[width=0.4\textwidth]{TestThermal/Vgw.pdf}
%\caption{Liquid and gas superficial velocities at different times along the left vertical left branches (b1) and the right vertical + horizontal branches (b2) for the chair shaped thermal well test case.}
%\label{fig_test2_Vlg}
%\end{center}
%\end{figure}

\subsubsection{Thermal test case with cross flow}\label{sec:crossflow}

\begin{figure}[H]
\begin{center}
\includegraphics[width=0.1\textwidth]{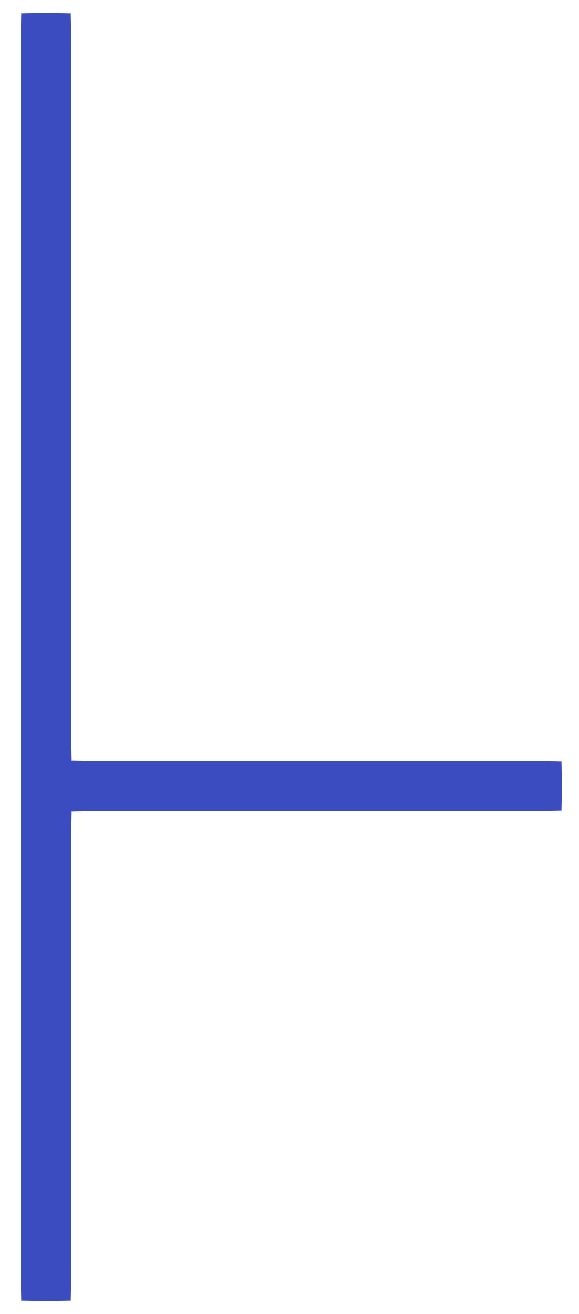}   
\caption{T shaped well with one junction and three branches of sizes 600 m (upper vertical branch), 400 m (lower left vertical branch), and 400 m (horizontal branch).}
\label{fig_test3_Twell}
\end{center}
\end{figure}

We consider the same liquid vapor thermal model as in the previous test case along the T shaped production well illustrated in Figure \ref{fig_test3_Twell}.
The DFM model \cite{shi2005} is again used with the same parameters as in the previous test case. In order to obtain cross flow, the reservoir properties at the two leaf nodes are set up as follows. The bottom left leaf will act as an injection point and is set up with the reservoir pressure $9~10^6$ Pa, the reservoir temperature $520$ K, the 
Darcy well index $WI^D = 10^{-12}$ m and the Fourier well index $WI^F = 100$ J.s$^{-1}$.K$^{-1}$. The right leaf will act as a feed zone and is set up with the reservoir pressure $7~10^6$ Pa, the reservoir temperature $500$ K, the 
Darcy well index $WI^D = 10^{-12}$ m and the Fourier well index $WI^F = 100$ J.s$^{-1}$.K$^{-1}$. 

The well is monitored as in the previous test case with the minimum well head pressure set to $\bar p_\omega = 5 ~10^5$ Pa and the maximum well mass flow rate set to $\bar q_\omega = 15$ Kg.s$^{-1}$. The initial temperature is set to $350$ K along the well and the initial pressure is hydrostatic at liquid state with pressure set to $5 ~10^5$ Pa at the head node.

The well is meshed using a uniform discretization of each of the three branches with $40$ edges for the bottom vertical branch, $40$ edges for the horizontal branch and $60$ edges for the top vertical branch. The simulation time is fixed to $t_F=2500$ s to reach the stationnary state and the time stepping is set up with an initial time step of $0.1$ s and a maximum time step of $40$ s. Using the same nonlinear solver setting as in the previous test case, the simulation runs in $147$ time steps with $5$ time step failures and a total number of $1056$ Newton iterations.

\begin{figure}[H]
\begin{center}
  \includegraphics[width=0.2\textwidth]{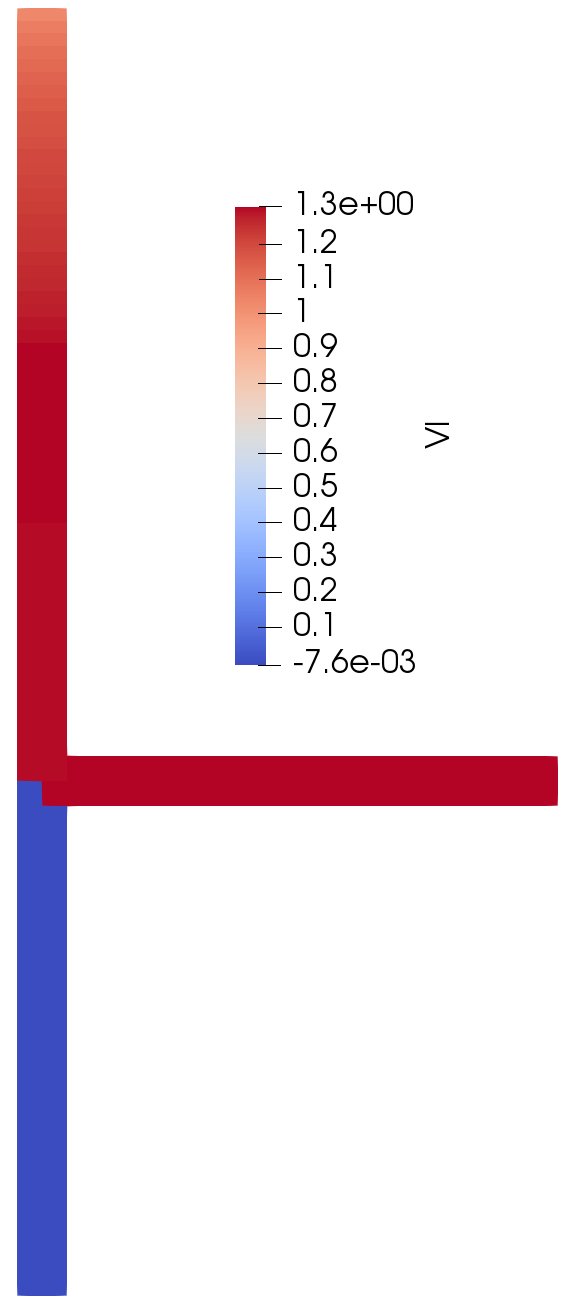}
  \includegraphics[width=0.2\textwidth]{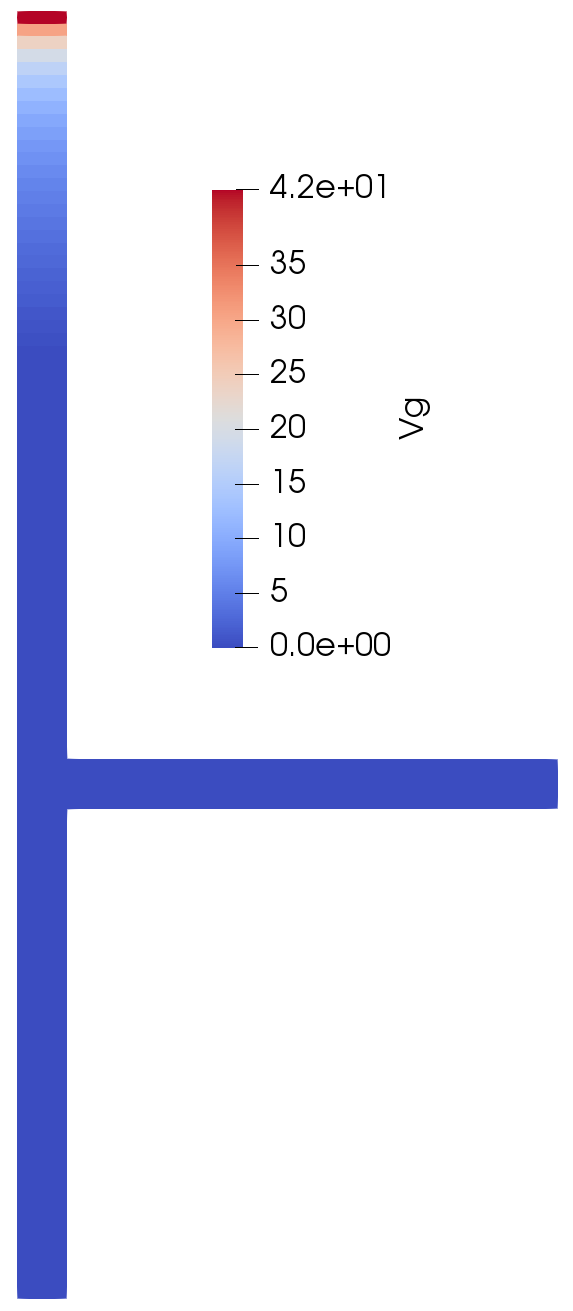}
\caption{Liquid and gas superficial velocities (m.s$^{-1}$) at final time for the cross flow thermal well test case.}
\label{fig_test3_Vlg_paraview}
\end{center}
\end{figure}

\begin{figure}[H]
\begin{center}
  \includegraphics[width=0.4\textwidth]{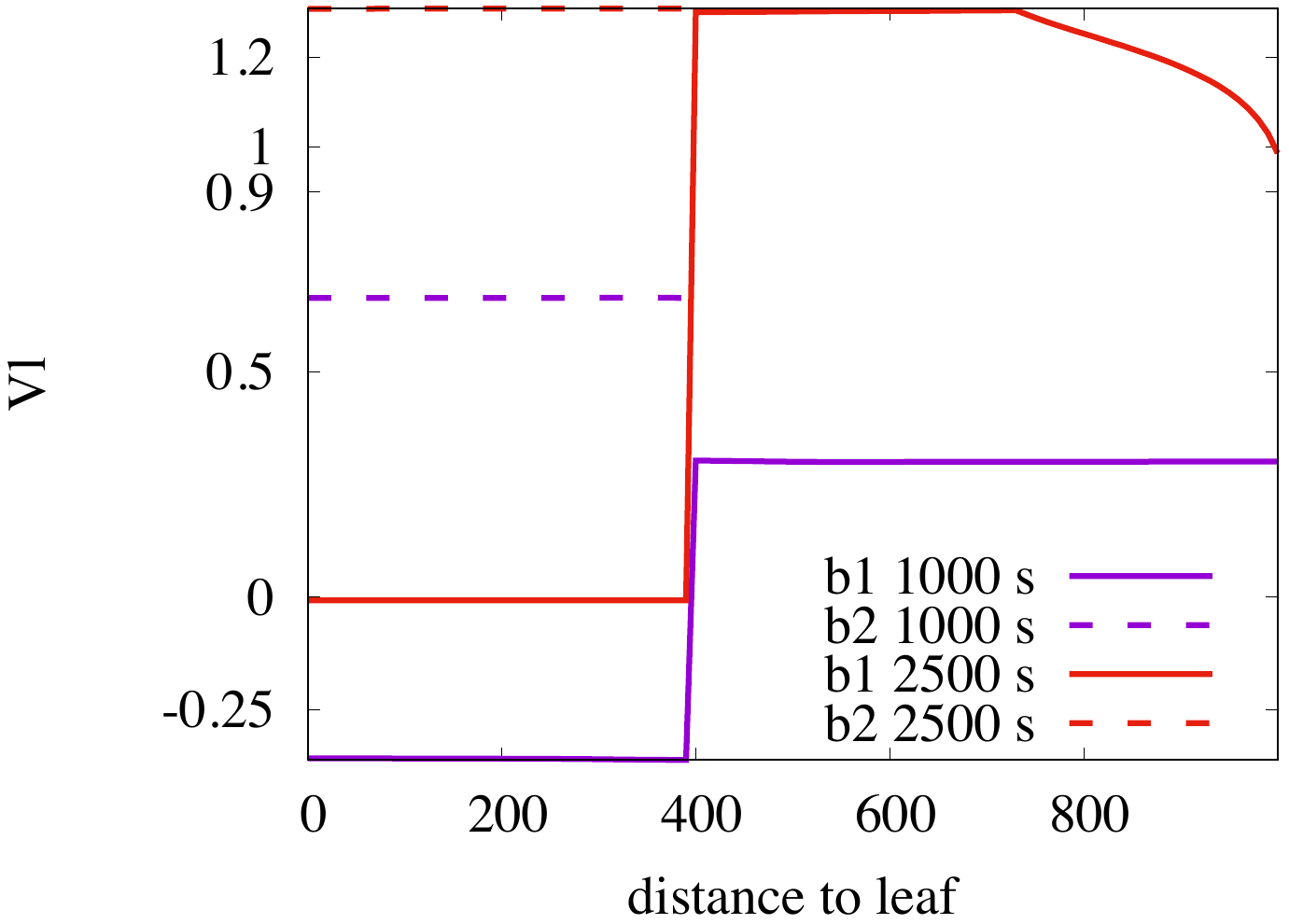}
  \includegraphics[width=0.4\textwidth]{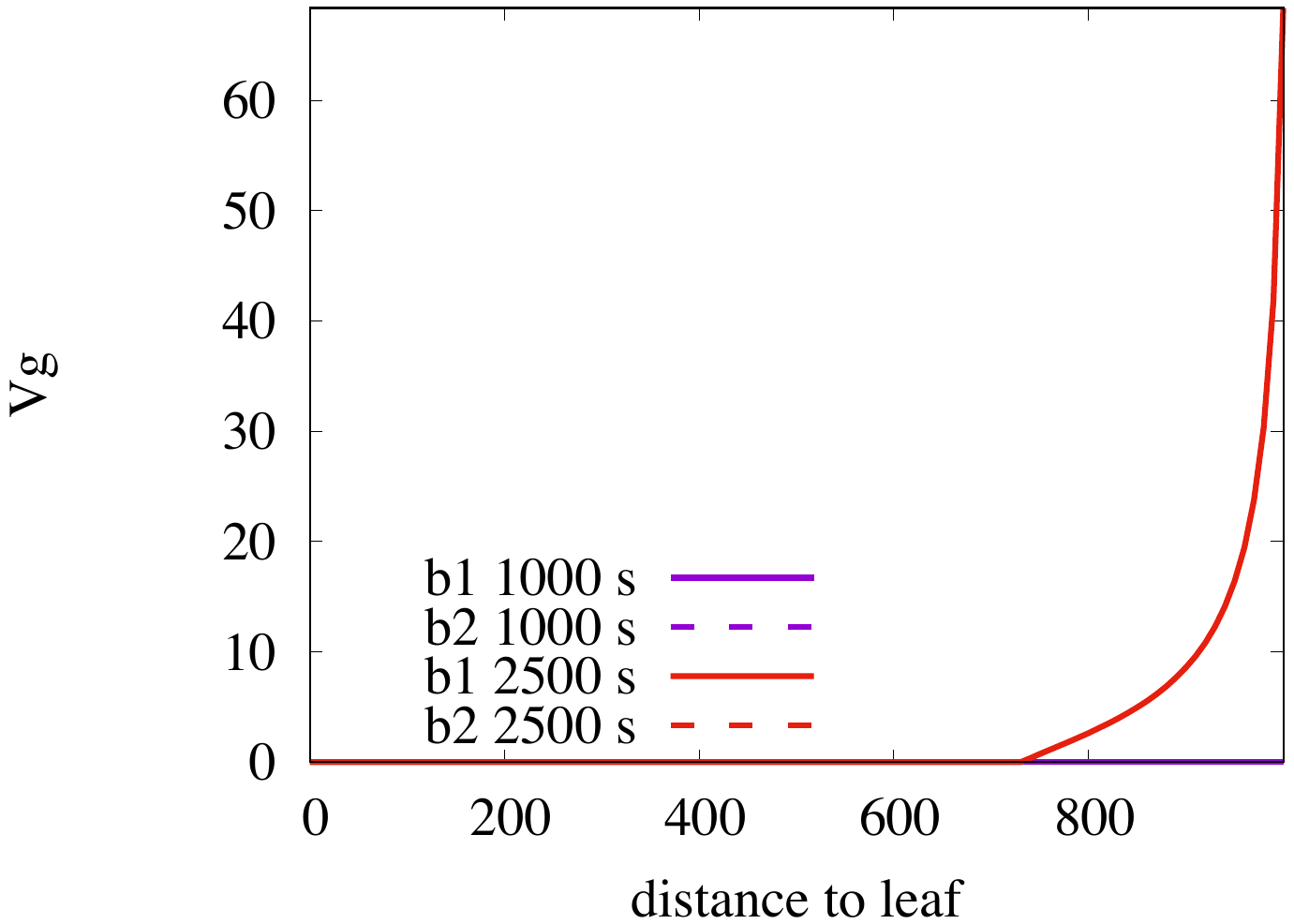}
\caption{Liquid and gas superficial velocities at different times along the left vertical left branches (b1) and the right horizontal branch (b2) for the cross flow thermal well test case.}
\label{fig_test3_Vlg}
\end{center}
\end{figure}

Figure \ref{fig_test3_Vlg_paraview} exhibits the final superficial liquid and gas velocities and Figure \ref{fig_test3_Vlg} plots the liquid and gas superficial velocities along the vertical (b1) and horizontal (b2) parts of the well  at times $1000$ s and $2500$ s. These figures clearly exhibit the cross flow between the right leaf acting as a feed zone and the bottom leaf acting as an injection point.

\begin{figure}[H]
\begin{center}
  \includegraphics[width=0.2\textwidth]{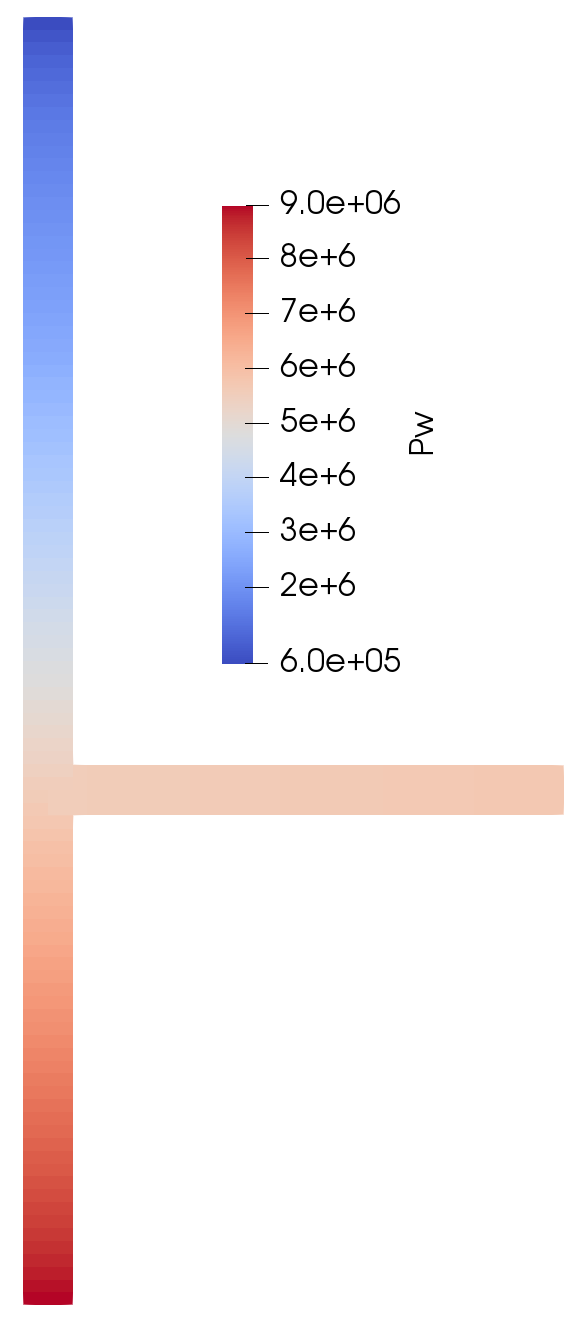}
  \includegraphics[width=0.2\textwidth]{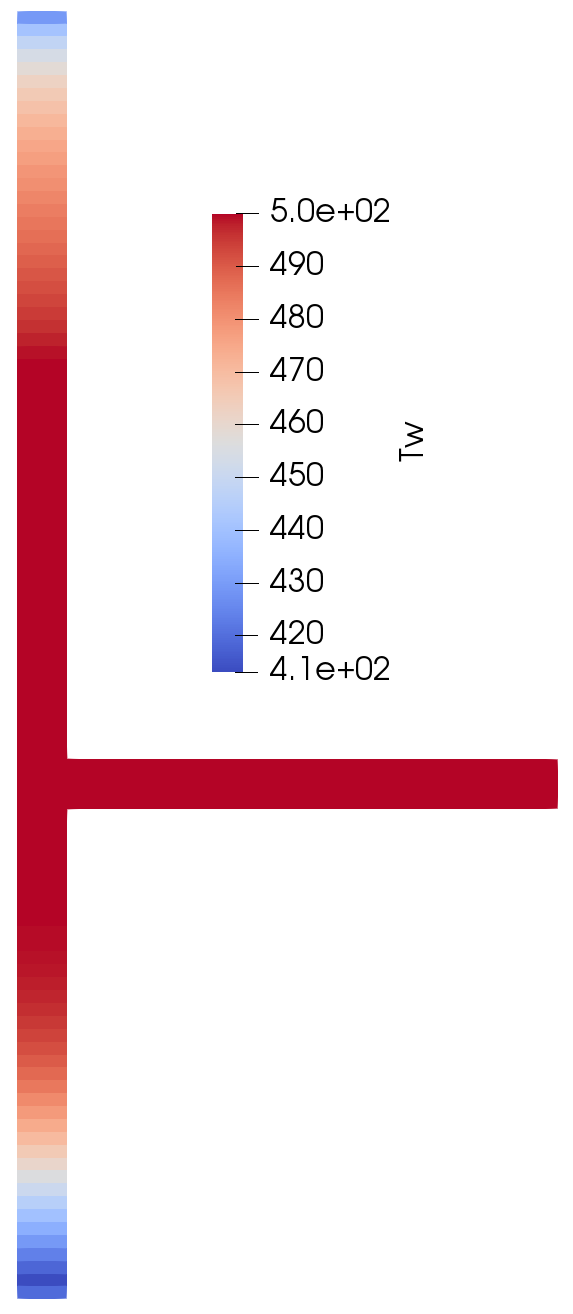}
  \includegraphics[width=0.2\textwidth]{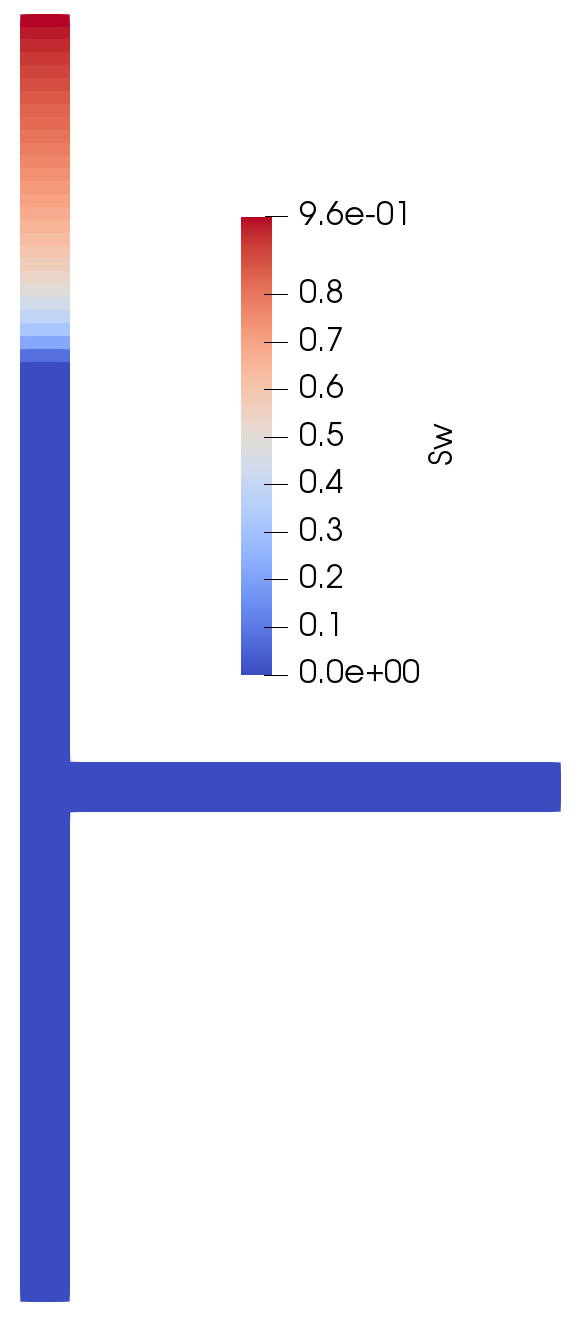}   
\caption{Pressure (Pa), temperature (K), gas saturation solutions at final time for the cross flow thermal well test case.}
\label{fig_test3_PTS_paraview}
\end{center}
\end{figure}

\begin{figure}[H]
\begin{center}
  \includegraphics[width=0.4\textwidth]{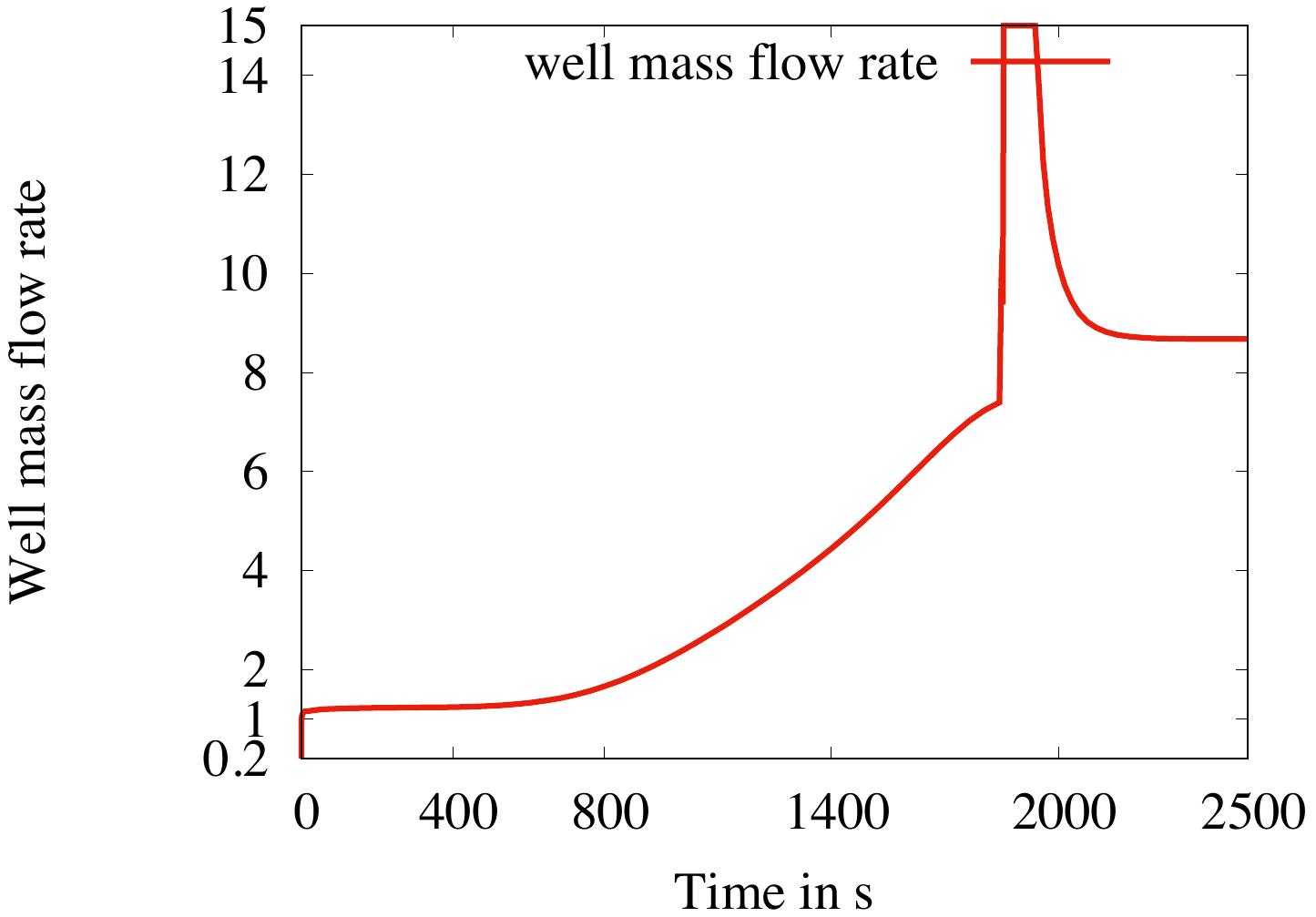}
  \includegraphics[width=0.4\textwidth]{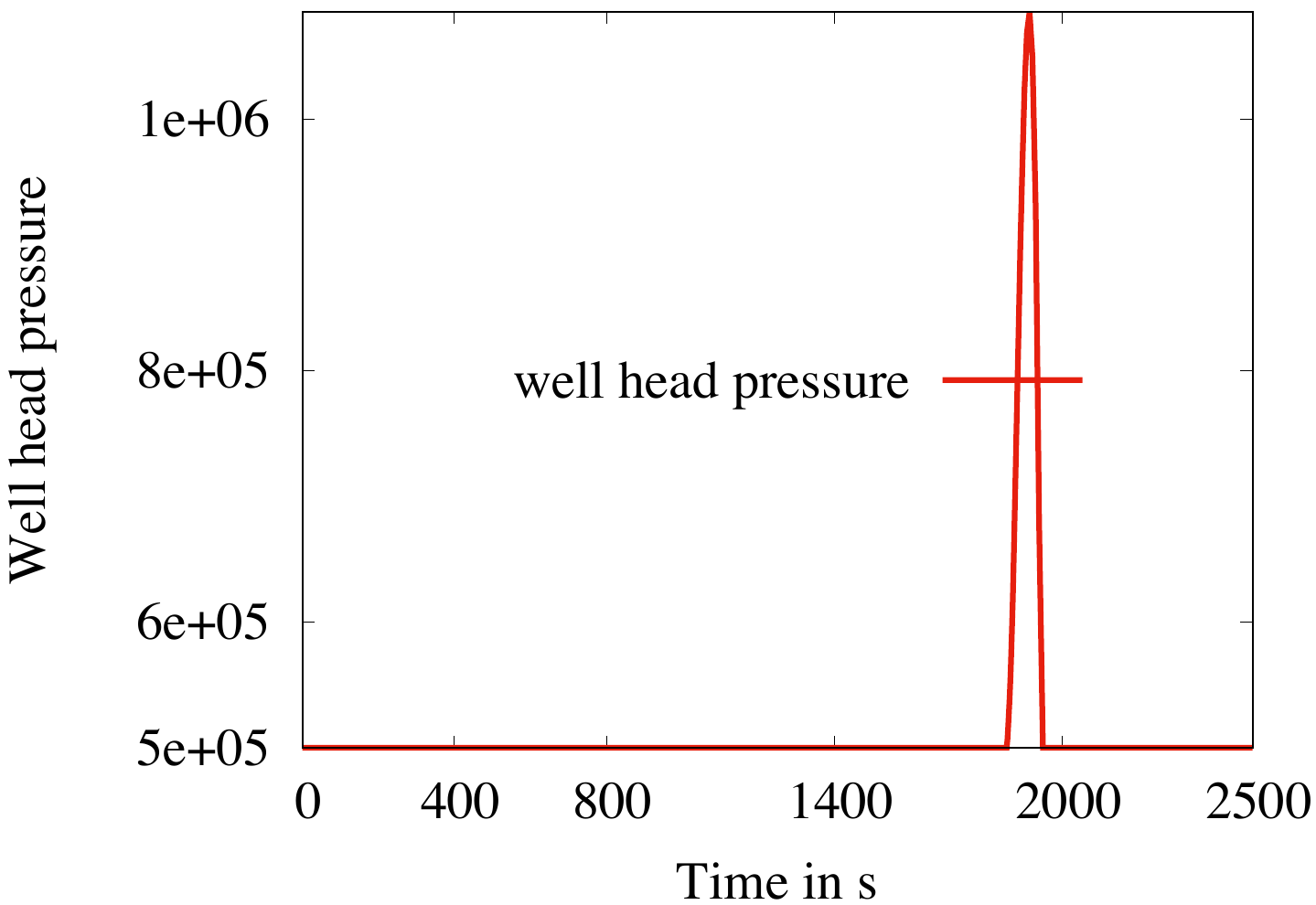}
  \includegraphics[width=0.4\textwidth]{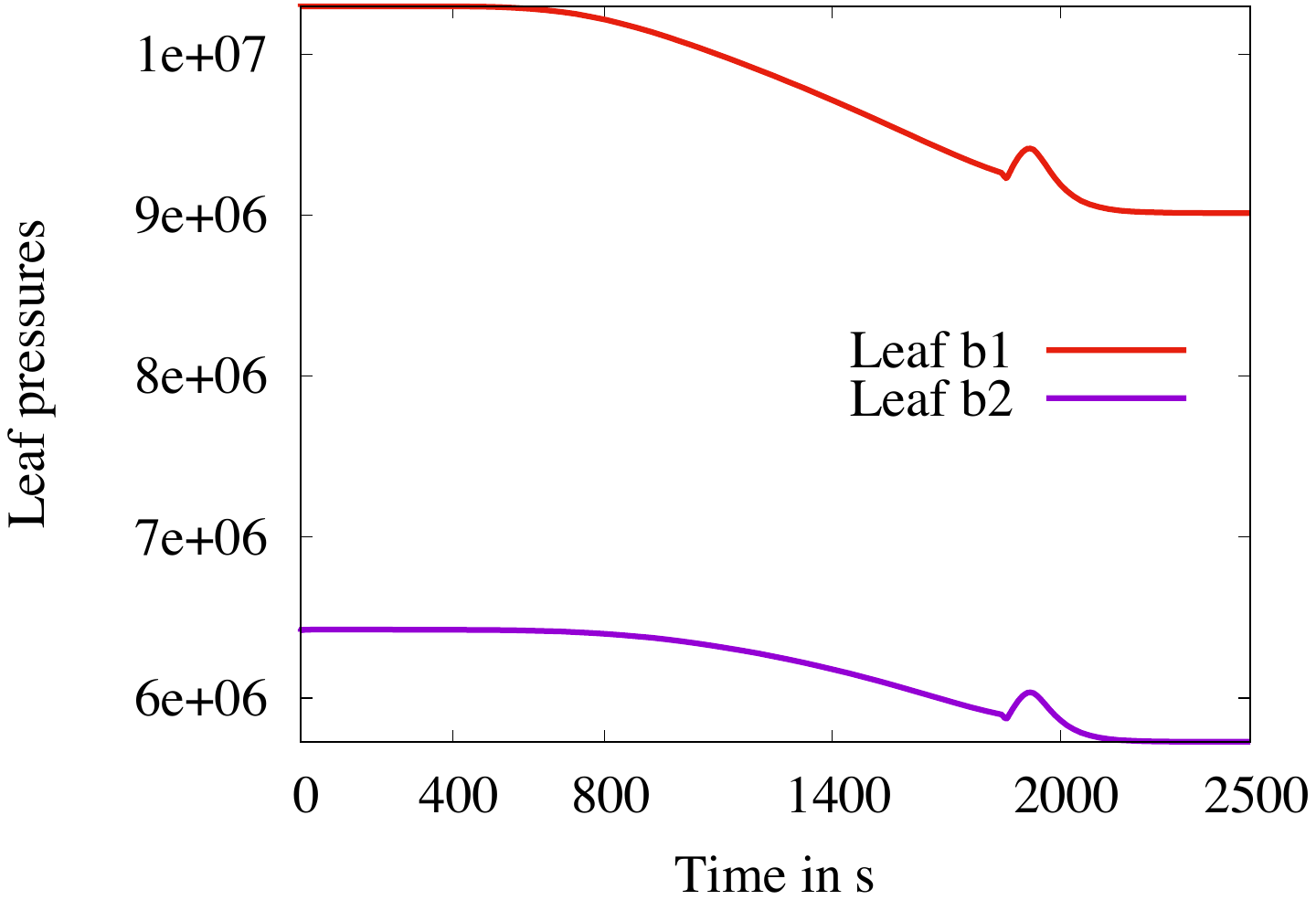}
   \includegraphics[width=0.4\textwidth]{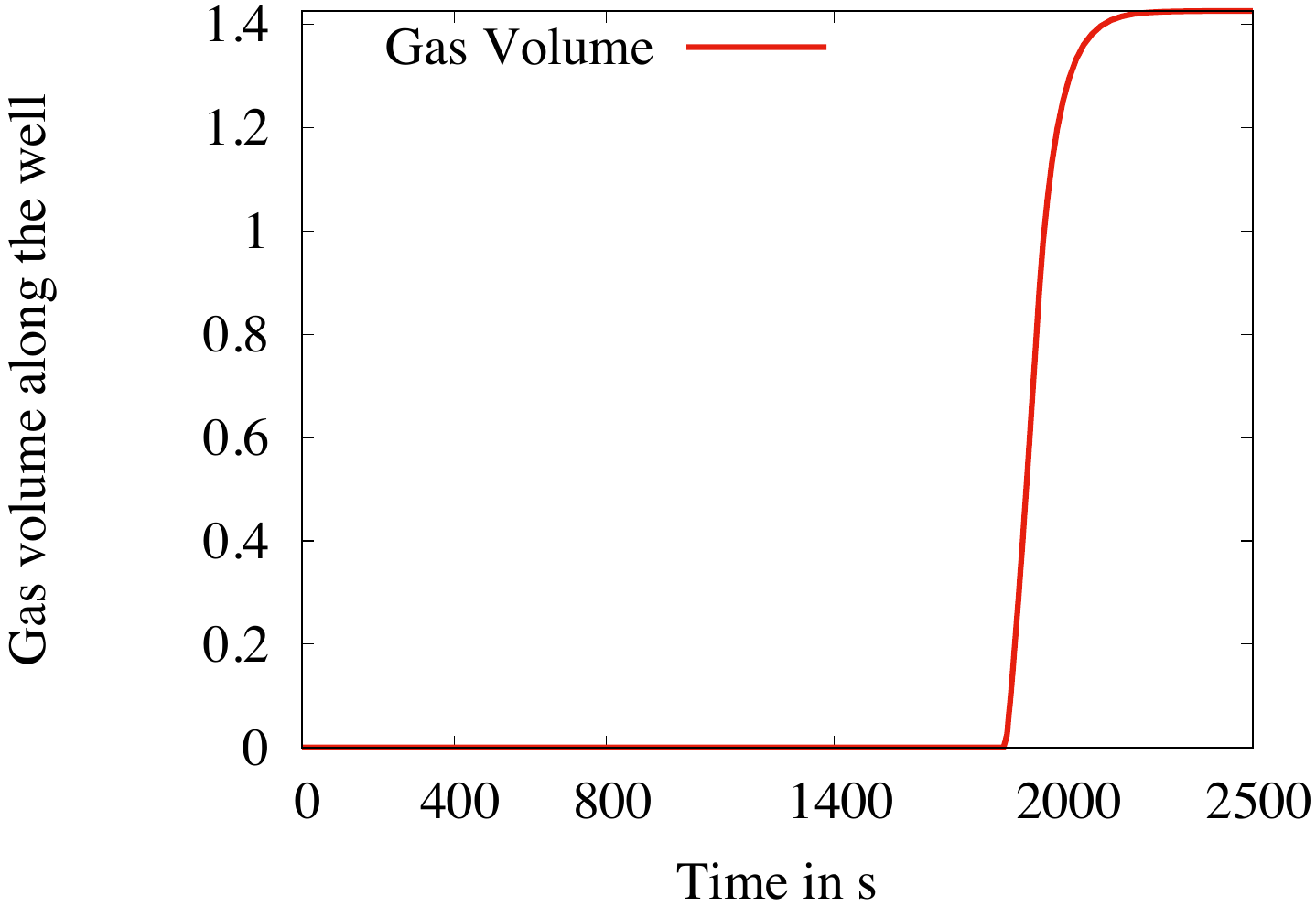}   
\caption{Well mass flow rate, well head pressure, leaf pressures, and gas volume in the well as a function of time for the cross flow thermal well test case.}
\label{fig_test3_histories}
\end{center}
\end{figure}

Figure \ref{fig_test3_PTS_paraview} exhibits the final pressure, temperature, and gas saturation. The high temperature front propagates from the right leaf to both the top and the bottom sides of the vertical part of the well. 
The rise of the hot temperature front at the low pressure top side of the well induces the appearance of the gas phase starting at roughly $t=1800$ s (see also Figure \ref{fig_test3_histories}). Figure \ref{fig_test3_histories} plots the time histories of the well pressure $\bar p_\omega$, the mass flow rate $\bar q_\omega$, the leaf pressures and the gas volume inside the well. It is shown that the well is monitored at the minimum well head pressure until the gas phase appears reducing the weight of the column which reduces the pressure at the right leaf node and allows to reach the maximum mass flow rate.
The well rapidly switches back to the minimum pressure monitoring due to the increase of the wall friction induced by the high gas velocity.

\begin{figure}[H]
\begin{center}
  \includegraphics[width=0.4\textwidth]{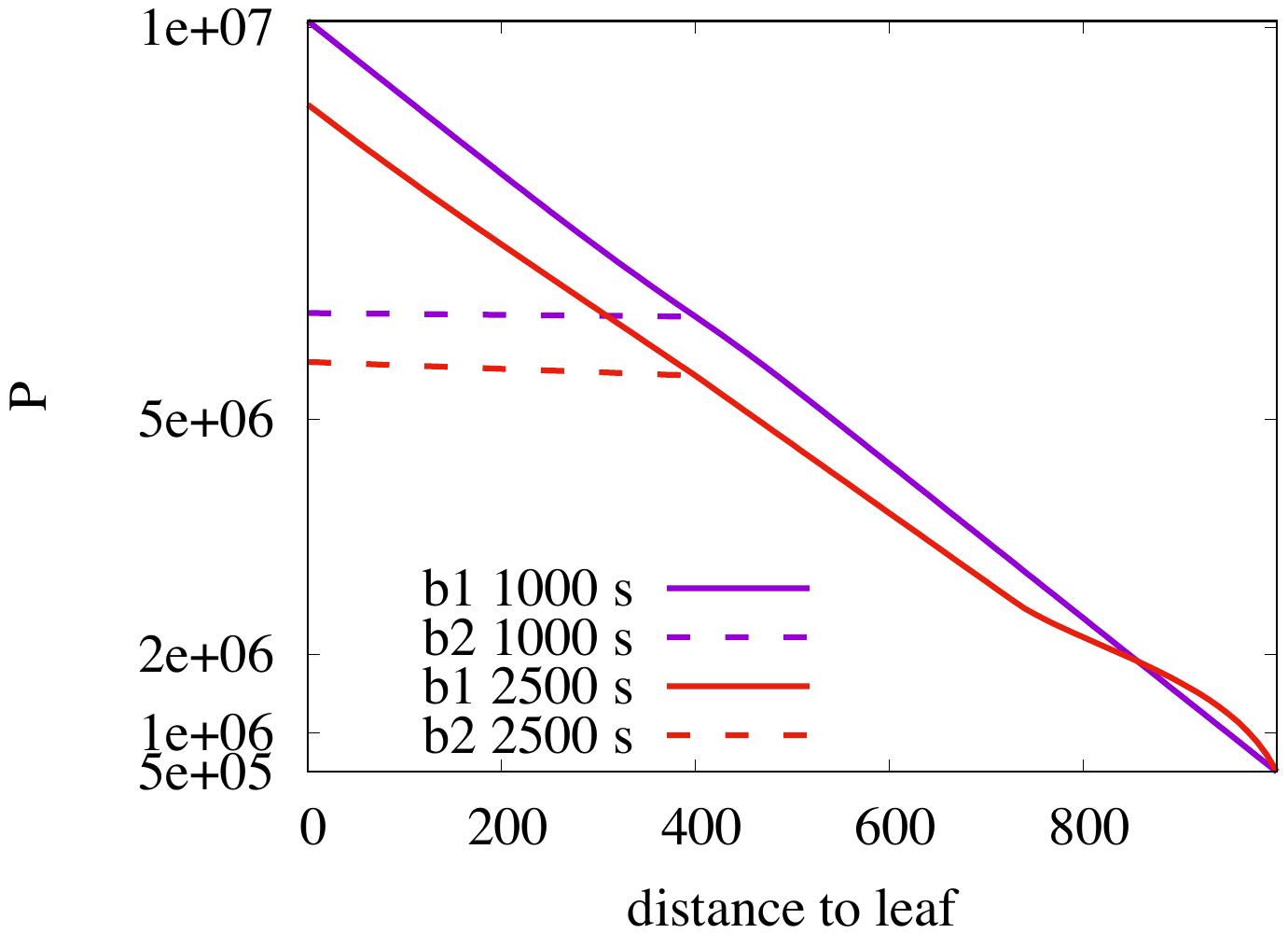}
  \includegraphics[width=0.4\textwidth]{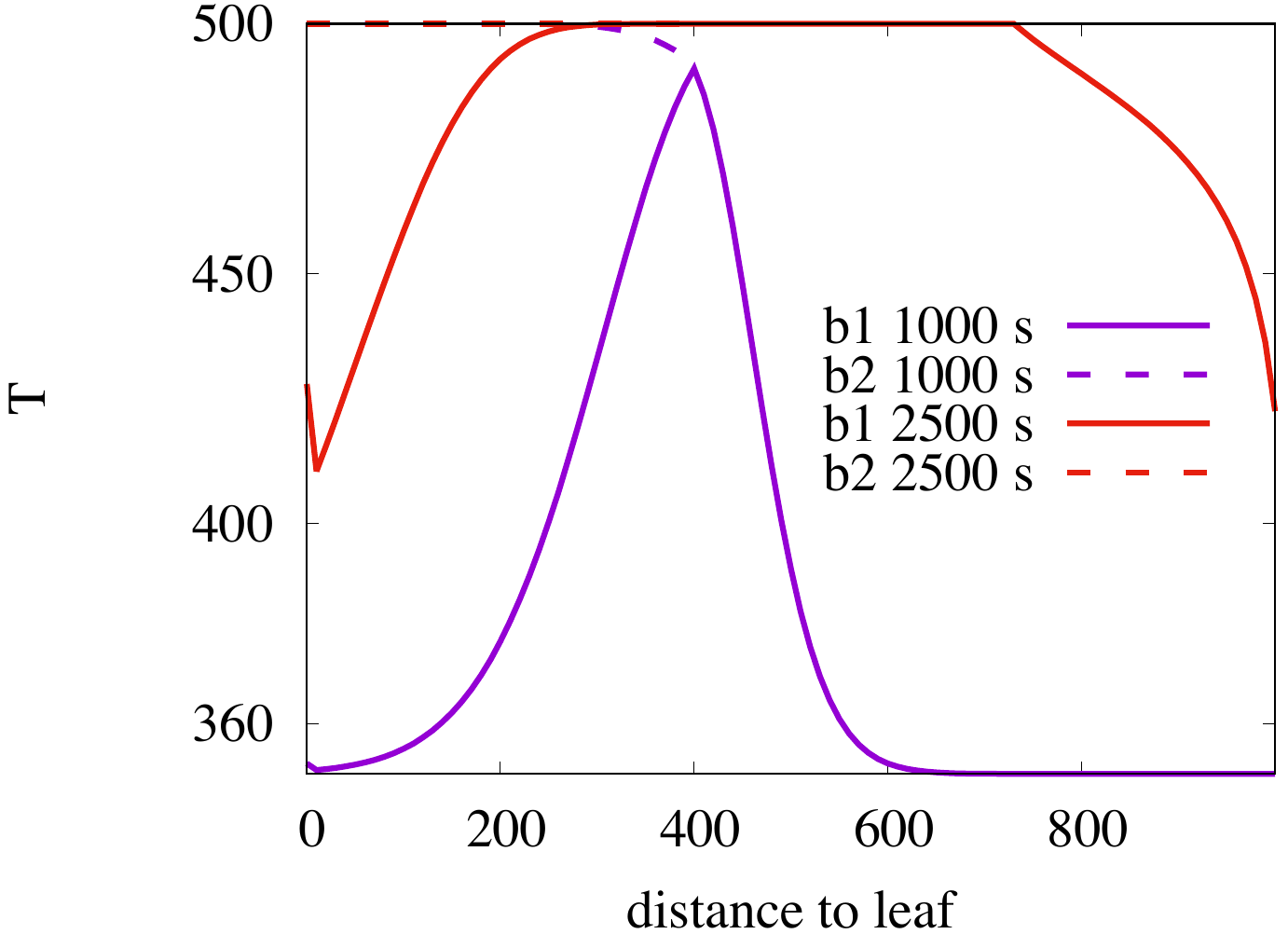}
  \includegraphics[width=0.4\textwidth]{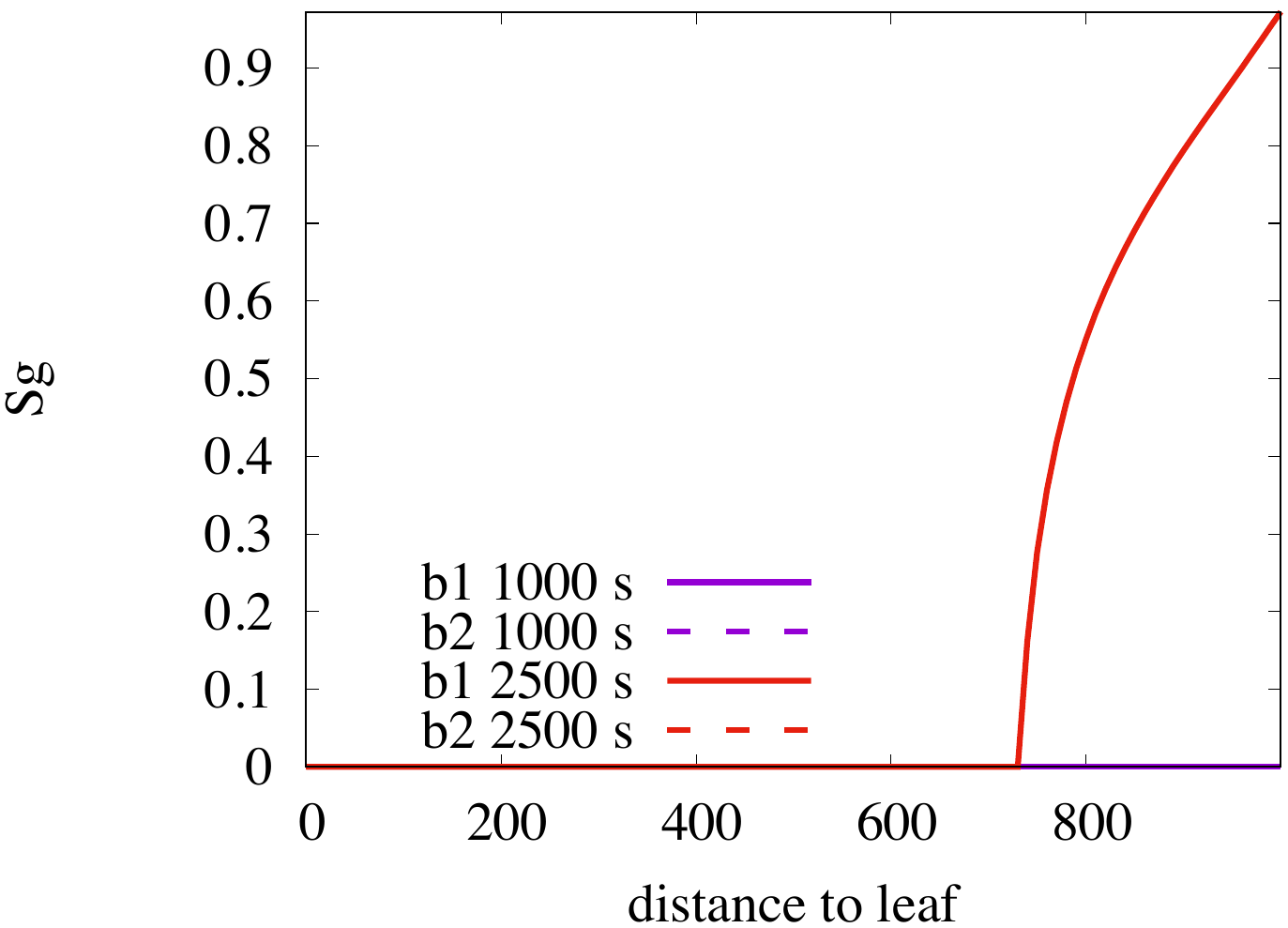}   
\caption{Pressure, temperature and gas saturation at different times along the left vertical left branches (b1) and the right horizontal branch (b2) for the cross flow thermal well test case.}
\label{fig_test3_PTS}
\end{center}
\end{figure}

Figure \ref{fig_test3_PTS}  exhibits the plots of the pressure, temperature, gas saturation along the vertical (b1) and horizontal (b2) parts of the well at times 1000 s and 2500 s showing the propagation of the temperature front on both sides of the vertical part of the well, as well as the decrease of the pressure along the well as a result of the temperature increase.

%\subsection{Validation of the fully coupled model}

\subsection{One production well coupled with a reservoir} \label{sec:fullycoupled}

The objective of this test case is to validate the coupled reservoir multi-segmented well model by comparison  with the results obtained using a simpler well model based on a single implicit unknown.  
We consider single component ${\rm H}_2{\rm O}$ two-phase (liquid and vapor), non-isothermal flow.
The DFM model \cite{shi2005} parameters  are the same as the ones used in the test case  \ref{sec:chairwell.test}.  
However, 
the Darcy-Forccheimer friction law \eqref{eq_darcy_forchheimer} is set up with the friction parameter 
$f_{q} \in \{ 0.001, 0.06\}$
where the case $f_{q}= f_{q_1} = 0.001$ corresponds to a very small friction coefficient.
 This choice will be useful for comparison with previous published results in the literature without friction as explained below. The case $f_q = f_{q_2} = 0.06$ corresponds a typical order of magnitude. 
 The geothermal reservoir is defined by the domain $\Omega = (-H, H)^2\times(0,H_z)$
where $H=1000$ m and $H_z=200$ m, and we consider one vertical producer well 
along the line $\{(x, y, z) \in \Omega \,|\, x = y = 0\}$ of radius $r_{\omega}=0.1$ m.
The reservoir is assumed homogeneous with isotropic permeability 
$\K = k {\mathbb I}, k = 5\times 10^{-14} \text{ m}^2$ and porosity $\phi = 0.15$.
It is assumed to be initially saturated with pure water in liquid phase. 
The internal energy, mass density and viscosity of water in the liquid and gas phases  are defined by analytical laws as functions of the pressure and temperature (refer \cite[Section 4.1]{ACJLM23} for details).
The vapour pressure $P_{sat}(T)$ is given in Pa by the Clausius-Clapeyron equation 
$$
p_{\rm sat}(T) = 100 \exp\left(46.784 - \frac{6435}{T} - 3.868 \,\, log(T)\right). 
$$
The reservoir thermal conductivity is fixed to $\lambda^r = 2 \;\text{W}.\text{m}^{-1}.\text{K}^{-1}$,
 and the rock volumetric heat capacity is given by $C_{s} = 1.6\; \text{MJ}.\text{K}^{-1}.\text{m}^{-3}$ with $E^s(p,T) = C_s T$.  
The relative permeabilities are set to
$k_{r}^\alpha(s^\alpha) = (s^\alpha)^2$ for both phases $\alpha\in \{\l,\g\}$. 
The gravity vector is as usual ${\bf g}=(0,0,-g_z)$ with $g_z=9.81 \text{ m}.\text{s}^{-2}$.
The simulation consists in two stages.
At the first one, the well is closed and  we impose a Dirichlet boundary condition at the top
of the domain prescribing the reservoir pressure and the temperature equal to $p^r=4$ MPa
and $T^r = (p_{\rm sat})^{-1}(p^r)- 1\;\text{K}$; respectively,
and homogeneous Neumann boundary conditions are set at the bottom and at the sides of the domain.
{The choice of the initial temperature $T^r$ just below the saturated vapor temperature is made in order to make the gas phase appear at the beginning of the production during the second stage.}
The first stage is run until the simulation reaches a stationary state with the liquid phase only, a constant temperature and an hydrostatic pressure depending only on the vertical coordinate. 

For the second stage, homogeneous Neumann boundary conditions are prescribed at the bottom and at the top of the domain $\Omega$,
but Dirichlet boundary conditions for the pressure and temperature are fixed at the sides of the domain to the ones at the end of stage one. 
The well is  set in an open state, i.e.,
it can produce,
and its monitoring conditions are defined by the minimum bottom hole pressure
$\bar{p}_\omega= 1$ bar (never reached in practice) and the maximum total mass flow rate $\bar{q}_\omega= 200  \text{ ton}.\text{hour}^{-1}$. The second stage is run on the time interval  $(0,t_F)$ with $t_F=30$ days. 
Moreover, at the top of the reservoir and at the root of the well, there is no coupling between them during both stages.
The entire simulation runs on uniform Cartesian mesh of size $n_x\times n_y \times n_z$ to discretize the domain $\Omega$ with 
$(n_x,n_y,n_z) = (20,20,10)$.
The well indexes $WI^D$ and $WI^F$ introduced in Appendix \ref{sec:peaceman} are computed at each node of the well following \cite{Xing2018}.

In order to validate the coupled model, we compare
the numerical results against the ones published in \cite[Section 4.1]{ACJLM23} using the same thermodynamic parameters, the same initialization and stages,  and the same mesh. The model presented in \cite{Xing2018} is based on the following assumptions on the well model:
\begin{itemize}
\item[(i)] the wall friction, the thermal conduction and the transient terms are all neglected all along the well,
 \item[(ii)] the thermal conduction is neglected between the reservoir and the well,  
\item[(iii)] there is no cross flow and both the gas and liquid velocities are oriented in the same direction all along the well,   
\item[(iv)] the pressure drop along the well is computed based on an explicit approximation of the mean density $\rho^{\rm m}$. 
\end{itemize}
In that case, the well model can be shown (see e.g. \cite{Xing2018}) to reduce to a single implicit unknown, the well head node pressure $p_\omega$, and to a single equation,  the monitoring conditions, fully coupled to the reservoir system. The computation of the mean density  $\rho^{\rm m}$ along the well is based on a gas liquid flash computation providing the well temperature and gas saturation at each well node. This flash computation is based on the lagged in time values of the well pressures and molar and energy flow rates $q^{r\rightarrow \omega}_{\s,i}$, $q^{r\rightarrow \omega}_{\s,\energy}$. It usually assumes a zero slip law $u^\g = s^\g u^{\rm m}$ as detailed in \cite{Xing2018}, but it can also easily account for a general slip law  \eqref{slip_law} by a simple modification of the gas liquid flash computation.     
In order to investigate further the effects of the DFM slip law, we consider in this test case the model proposed
in \cite{Xing2018} with both a zero slip law and its enhancement using the same slip law as the one used in the multi-segmented well model. For the sake of completeness, the modified flash computations implied by the non zero slip law is detailed in Appendix \ref{sec:siuwell}.

\begin{figure}[H]
 \begin{center}
   \includegraphics[width=0.45\textwidth]{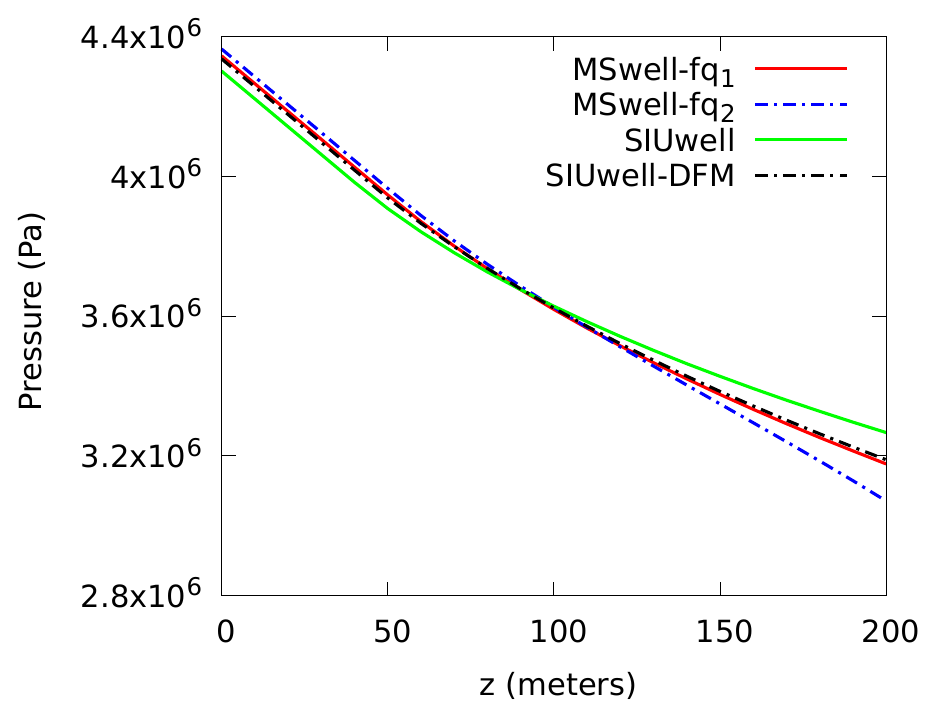}
   \includegraphics[width=0.45\textwidth]{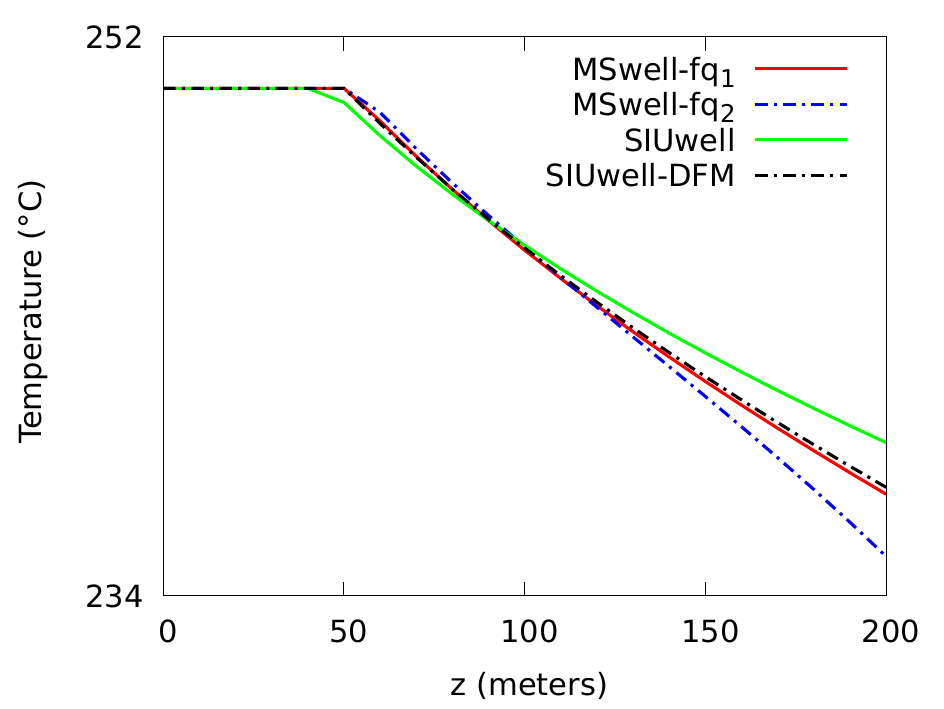}\\
   \includegraphics[width=0.45\textwidth]{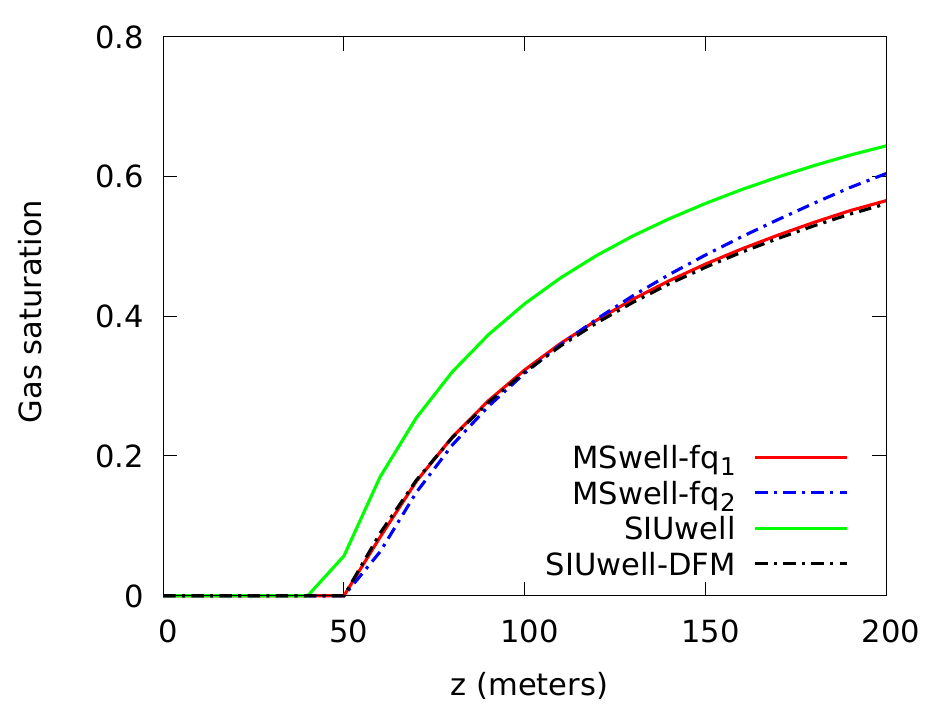}
     \caption{Pressure in Pa,  temperature in $^\circ$C  and gas saturation along the well at final time using the three different  models.}
     \label{fig_pTs_along_well}
  \end{center}
 \end{figure}

\begin{figure}[H]
\begin{center}
  \includegraphics[width=0.45\textwidth]{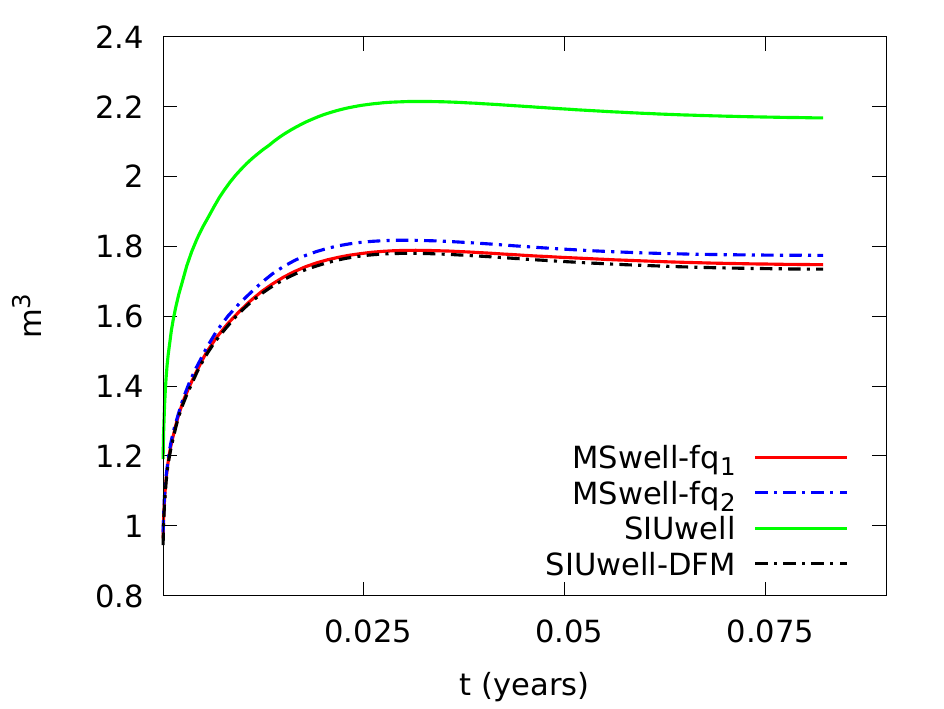}
  \includegraphics[width=0.45\textwidth]{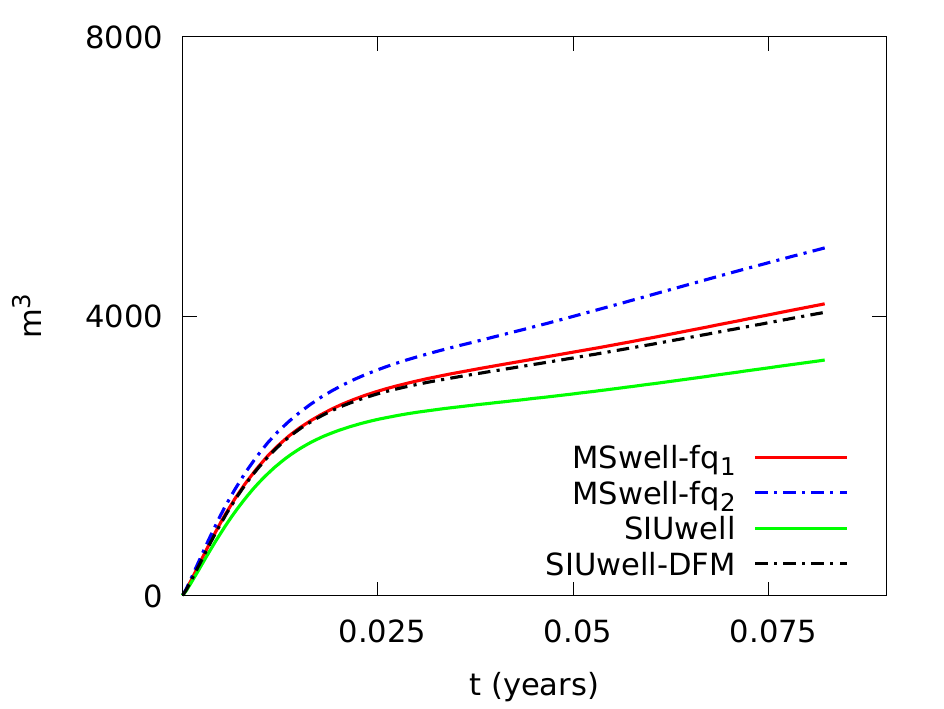}
    \caption{Total gas volume inside the well (left) and inside the reservoir (right) as a function of time  using the three different models.}
    \label{fig_gas_vol_well_res}
\end{center}
\end{figure}

To  make easier the comparison of the results obtained by each model, we label by \texttt{MSwell-fq}$_1$ and \texttt{MSwell-fq}$_2$
the results obtained by the proposed Multi-Segmented well model \eqref{eq_well_discrete} with $f_{q_1}=0.001$ and $f_{q_2}=0.06$, respectively. 
We label by \texttt{SIUwell} the ones obtained by
the Single Implicit Unknown well model proposed in \cite{Xing2018},
and we label by \texttt{SIUwell-DFM} its enhancement including the DFM  slip law \eqref{slip_law} in the gas liquid flash computation. Figure \ref{fig_pTs_along_well} 
 compares the pressure, the temperature and the gas saturation along the well; respectively, at final time $t_F$ using this three different models.
Figure \ref{fig_gas_vol_well_res}  show the total 
volume of gas inside  the well and  the reservoir as functions of time for the three models.
It can be noticed from those figures that the results obtained with the
MSwell model \eqref{eq_well_discrete} with $f_{q_1}=0.001$ are very close to the ones obtained  using the well model of \cite{Xing2018} provided that the DFM slip law is taken into account. This is expected for such configuration with no cross flow, no significant thermal conduction losses between the well and the reservoir, no significant wall friction, and no significant transient effects at the reservoir time scale. 
On the other hand, the more realistic value $f_{q_2}$ of the friction coefficient induces a larger pressure drop which results in significant variations between the  \texttt{MSwell-fq}$_2$ and the \texttt{SIUwell-DFM} models.  

Table \ref{tab_numbeh} shows the numerical efficiency of the nonlinear solver for all  models using the same mesh for the second stage of the simulation. We denote by $N_{\Delta t}$  the number of time steps and by $N_{\text{Newton}}$ the average number of Newton iterations per time step. It exhibits as expected the  higher number of Newton iterations for the multi-segmented well models compared with the single implicit unknown well models, as a result of stronger nonlinearities. We note also that the nonlinear convergence of the multi-segmented well model is sensitive to very small (not physical) friction coefficients due to an increased stiffness of the mixture velocity as a function of $\Delta_\welledge \Phi$ in such cases (see \eqref{eq_um}). 

\begin{table}[H]
\begin{center}
\begin{tabular}{|c|c|c|c|c|} \hline
  Model  & $\# \cells $  & $N_{\Delta t}$ &   $N_{\text{Newton}}$   \\ \hline
  \texttt{MSwell-fq}$_1$ & 32000 & 121   & 10.12 \\
  \texttt{MSwell-fq}$_2$ & 32000 & 120   & 6.15 \\
  \texttt{SIUwell}  & 32000 & 133   & 1.56 \\
  \texttt{SIUwell-DFM}   & 32000 & 133   & 1.58 \\\hline
\end{tabular}
\caption{
Numerical behavior of the second stage of the simulation for different models using the same mesh.
$N_{\Delta t}$ is the number of time steps, and $N_{\text{Newton}}$ the average number of Newton iterations per time step.}
  \label{tab_numbeh}
\end{center}
\end{table}

\section{Conclusion}
A  numerical model for two-phase compositional non-isothermal flow in geothermal multi-branch wells is developed to simulate the flow and transport along the wells in geothermal field operations. The model combines the flexible thermodynamical Coats' formulation with the a drift flux hydrodynamical model leading to a system coupling the conservation equations for each component, momentum and energy.

The discretization is fully implicit in time and based on a staggered finite volume scheme in space, with node centred control volumes for components and energy conservations and edge control volumes for the momentum conservation. The numerical fluxes combine a monotone flux approximation for the phase superficial velocities on the full range of gas saturation with an upwind approximation of the phase molar fractions, density and enthalpy. The interplay between these two parts of the fluxes is a key ingredient of the stability of the scheme. 

The nonlinear solver benefits from the elimination of the well superficial velocities and flow rates in the Newton linearization process leading to a Jacobian system with the same primary unknowns both for the well and the reservoir, which makes the coupling between both models much easier. 

This numerically robust well model can be used to simulate complex transient flows such as those occurring during the start-up of a geothermal well, as well as cross-flow configurations along the wells. 
It can be used either by modeling the reservoir with source terms using well indexes, or fully  coupled with a subsurface flow model. Different configurations (transient, stationary, coupled or not with a reservoir flow model) have been tested to validate the model performance.

Future work will focus on the implementation of alternative drift flux models in particular to take downward flows into account, and on coupled well-reservoir modeling for large scale industrial case studies.

\section{Appendices} \label{sec:appendices}

\subsection{Appendix I: example of monotone two-point flux for the DFM model from \cite{shi2005}}
\label{sec_monotoneflux}
The gas liquid DFM model introduced in \cite{shi2005} is based on the following choices of the drift velocity and profile parameter.

Let $\sigma_{\g\l}$ be the gas liquid interfacial tension, and let us define the following characteristic velocity for the rise of a gas bubble in a liquid: 
$$
U_c = \( \sigma_{\g\l} g {(\rho^\l - \rho^\g) \over (\rho^\l)^2 } \)^{1\over 4}. 
$$

The profile parameter is defined by
\begin{equation}
\label{eq_C0}
\left\{\begin{array}{r@{\,\,}c@{\,\,}l}
& \dsp C_0(s^\g) = {A \over 1 + (A-1)\gamma^2 }, \\[2ex]
& \dsp \gamma = P_{[0,1]}\({\beta - B \over 1-B}\), \\[2ex]
& \dsp \beta = \max\(s^\g, F_\nu s^\g {|u^{\rm m}| \over V_{sgf}}\),\\[2ex]
& \dsp V_{sgf} = K_u \sqrt{\rho^\l \over \rho^\g}  U_c,
\end{array}\right.
\end{equation}
with constant parameters $A$, $B$ which must be such that $B < (2-A)/A$ to ensure that $s^\g C_0(s^\g) \leq 1$ and
$s^\g C_0(s^\g)$ non decreasing. Typically we set $A=1.2$, $B=0.3$. The Critical Kutateladze number $K_u \in [0,3.5]$ depends on the dimensionless diameter of the pipe defined by
$$
\hat D = \(  g {(\rho^\l - \rho^\g) \over \sigma_{\g\l} } \)^{1\over 2} D, 
$$
where $D$ is the diameter of the pipe. 
The constant parameter $F_\nu$ is typically set to $1$ and the projection $P_{[0,1]}$ is defined by
$$
P_{[0,1]}(x) = 
\left\{\begin{array}{r@{\,\,}c@{\,\,}ll}
& 0 &\mbox{ if } x \leq 0,\\
& 1 &\mbox{ if } x \geq 1,\\
& x &\mbox{ if } x\in (0,1).  
\end{array}\right.
$$
The drift velocity times $s^\g$ is defined by 
\begin{equation}
\label{eq_Ud}
\left.\begin{array}{r@{\,\,}c@{\,\,}l}
s^\g U_{\rm d}(s^\g) = \dsp  G(s^\g) {\wt K}(s^\g)~U_c, 
\end{array}\right.
\end{equation}
with
$$
G(s^\g) = { (1-s^\g C_0(s^\g))  \over \dsp s^\g C_0(s^\g) \sqrt{\rho^\g\over \rho^\l}~+1 - s^\g C_0(s^\g)},
$$
and 
$$
{\wt K}(s^\g) = 
\left\{\begin{array}{r@{\,\,}c@{\,\,}ll}
& 1.53 ~s^\g &\mbox{ if } s^\g \leq a_1,\\[2ex]
& K_u C_0(s^\g) s^\g &\mbox{ if } s^\g \geq a_2,\\[2ex]
& \dsp 1.53 ~a_1 {s^\g-a_2 \over a_1-a_2 } + K_u C_0(a_2) a_2  {s^\g-a_1 \over a_2-a_1 }  &\mbox{ if } s^\g\in (a_1,a_2).  
\end{array}\right.
$$
Note that the linear interpolation between $a_1$ and $a_2$ is done in \cite{shi2005} on ${{\wt K}(s^\g)\over s^\g C_0(s^\g)}$ rather than on ${\wt K}(s^\g)$. It has been modified here to simplify the design of a monotone flux. 
The constant parameters $0 < a_1 < a_2 < 1$ are typically set to $a_1=0.2$ and $a_2=0.4$. 
We can check that $G(s^\g)$ is a non increasing function w.r.t. $s^\g$. We also assume in the following that
${\wt K}(s^\g)$ is a non-decreasing function w.r.t. $s^\g$ which is the case provided that the condition 
$$
1.53 ~a_1 \leq a_2 K_u C_0(a_2),
$$
is satisfied. In practice, it seems that this condition is not very restrictive. This condition guarantees that the following flux function 
$F^\g_\welledge$ is a monotone two-point flux: 

\begin{equation}
\label{eq_Fg}
\left\{\begin{array}{r@{\,\,}c@{\,\,}ll}
 F_\welledge^\g(u,v) &=& u ~C_0(u) (u^{\rm m}_\welledge)^+ + v ~C_0(v) (u^{\rm m}_\welledge)^- \\[2ex]
 && + G(v)  ~{\wt K}(u) (U_c o_\welledge)^+   + G(u) ~ {\wt K}(v) (U_c o_\welledge)^-. 
\end{array}\right.
\end{equation}

If $1.53 ~a_1 > a_2 K_u C_0(a_2)$, a Godunov scheme could still be computed for ${\wt K}(s^\g)$ since it is in that case a piecewise monotone function non decreasing on $(0,a_1)\cup(a_2,1)$ and non increasing on $(a_1,a_2)$.

\subsection{Appendix II: reservoir source terms}\label{sec:peaceman}

It is assumed that the radius of the wells are small compared to the cell sizes in the neighborhood of the well. It results that the Darcy flux between the reservoir and the well at a
given well node $\s\in \mathcal{V}_\omega$ is obtained using the Two Point Flux Approximation
$$
V_\s^{\alpha} = \WI^D_\s (p^{r,\alpha}_\s - p_\s),  
$$
where $p^{r,\alpha}_\s$ is the reservoir phase pressure at node $\s$. 
Fourier fluxes between the reservoir and the well are discretized in the same way using the Two Point Flux Approximation
$$
F_\s = \WI^F_\s (T^r_\s - T_\s),  
$$
with $T^r_\s$ denoting the reservoir temperature at node $\s$. The Well Indexes $\WI^D_{\s}$ and $\WI^F_{\s}$ are typically computed using Peaceman's approach (see \cite{Peaceman78,Peaceman83,CZ09,Xing2018}) and take into account the unresolved singularity of respectively the pressure and temperature solutions in the neighborhood of the well. Let us denote by $k^\alpha_{r,\s}(s^\alpha)$ the phase relative permeability at node $\s$ as a function of the phase saturation $s^\alpha$, by $s^{r,\alpha}_\s$ the reservoir saturation of phase $\alpha$, and by $c^{r,\alpha}_\s$ the reservoir molar fractions of phase $\alpha$.  

For any $a\in \R$, let us define $a^+ = \max(a,0)$ and $a^- = \min(a,0)$. 
The molar flow rates between the reservoir and the well $\omega$ at a given node $\s \in \mathcal{V}_\omega$ are defined by the following phase based
upwind approximation of the mobilities:
\begin{equation}
  \label{eq_q1sw}
  \begin{array}{r@{\,\,}c@{\,\,}l}  
    q^{r\rightarrow \omega}_{\s,\alpha,i} &=&  \dsp c_{\s,i}^{\alpha}{\zeta^\alpha_\s \over \mu^\alpha_\s} k^\alpha_{r,\s}(s_{\s}^\alpha)(V_\s^{\alpha})^- + c_{\s,i}^{\alpha}
      {\zeta^\alpha(p_\s^r,T^r_\s,c_\s^{r,\alpha}) \over \mu^\alpha(p^r_\s,T^r_\s,c_\s^{r,\alpha})} k^\alpha_{r,\s}(s_{\s}^{r,\alpha}) (V_\s^{\alpha})^+,\\[2ex]
    q^{r\rightarrow \omega}_{\s,i} &=& \dsp \sum_{\alpha\in \mathcal{P}_i} q^{r\rightarrow \omega}_{\s,\alpha,i}, 
    \end{array}
\end{equation}
and the energy flow rate is defined similarly by 
\begin{equation}
\label{eq_q2sw}
q^{r\rightarrow \omega}_{\s,e} = \sum_{\alpha\in \mathcal{P}} \( h^\alpha_\s (q^{r\rightarrow \omega}_{\s,\alpha})^- +
h^\alpha(p_\s^r,T^r_\s,c_\s^{r,\alpha}) (q^{r\rightarrow \omega}_{\s,\alpha})^+\)   + F_\s, 
\end{equation}
with $q^{r\rightarrow \omega}_{\s,\alpha} = \dsp \sum_{i \in \mathcal{C}^\alpha} q^{r\rightarrow \omega}_{\s,\alpha,i}$.

\subsection{Appendix III: computations of the well temperatures and saturations for the single implicit unknown well model}\label{sec:siuwell}

We detail in this Appendix the computations of the production well temperatures and saturations for the single implicit unknown well model of Subsection \ref{sec:fullycoupled} at given well pressures $p_\s$, molar $q^{r\rightarrow \omega}_{\s,\mass}$ and energy $q^{r\rightarrow \omega}_{\s,e}$ flow rates for $\s \in \mathcal{V}_\omega$. 
Using the assumptions (i) and (ii) stated in Subsection \ref{sec:fullycoupled} for the \texttt{SIUwell-DFM} well model, and considering the case of a single component ${\rm H}_2{\rm O}$, the conservation equations \eqref{cons_molar_discrete}-\eqref{cons_energy_discrete} reduce to 
$$
\begin{aligned}
 & \dsp  \sum_{\welledge\in \edges_\s^\omega}\sum_{\alpha\in \mathcal{P}}  - \kappa_{\welledge,\s} |S^\omega_\welledge| \zeta^\alpha_\welledge u^\alpha_\welledge = \dsp q^{r\rightarrow \omega}_{\s,\mass} - \delta_\s^{\s_\omega}  \sum_{\alpha\in \mathcal{P}}  |S^\omega_{\s_\omega}|  \zeta^\alpha_\s u^\alpha_\omega,  \\
& \sum_{\welledge\in \edges_\s^\omega}  \sum_{\alpha\in \mathcal{P}  } - \kappa_{\welledge,\s}  |S^\omega_\welledge| h^\alpha_\welledge \zeta^\alpha_\welledge u^\alpha_\welledge  
     \dsp  = \dsp q^{r\rightarrow \omega}_{\s,e} - \delta_\s^{\s_\omega} \sum_{\alpha\in \mathcal{P}} |S^\omega_{\s_\omega}| h^\alpha_\s \zeta^\alpha_\s u^\alpha_\omega,
\end{aligned}
$$
for $\s \in \mathcal{V}_\omega$. 
Summing these molar and energy conservation equations over all nodes $\s'' \underset{\omega}{\geq} \s$, and using assumption (iii) for the enthalpy and molar density upwind values, we obtain for all $\welledge=\s'\s\in \edges_\omega$ that
$$
\begin{aligned}
 & \dsp  \sum_{\alpha\in \mathcal{P}}  |S^\omega_\welledge| \zeta^\alpha(p_\s,T_\s) u^\alpha_\welledge =  \sum_{\s''\in {\cal V}_\omega | \s'' \underset{\omega}{\geq} \s} q_{\s'',\mass}^{r\rightarrow \omega} := Q_{\s,\mass}^\omega,   \\
&  \dsp \sum_{\alpha\in \mathcal{P}  }   |S^\omega_\welledge| h^\alpha(p_\s,T_\s) \zeta^\alpha(p_\s,T_\s) u^\alpha_\welledge  =  \sum_{\s''\in {\cal V}_\omega | \s'' \underset{\omega}{\geq} \s}  q^{r\rightarrow \omega}_{\s'',e} := Q_{\s,e}^\omega.   
\end{aligned}
$$
It results that the thermodynamical equilibrium at fixed well pressure $p_\s$, molar $Q_{\s,\mass}^\omega$ and energy $Q_{\s,e}^\omega$ provides the well temperature $T_\s$ and the well saturations $s_{\s}^{\alpha}$ at node $\s$ as follows.
Let us define the phase molar fractions $c_\s^\alpha$, $\alpha \in \mathcal{P}$ such that
$$
c^\alpha_\s  Q_{\s,\mass}^\omega =  |S^\omega_\welledge| \zeta^\alpha(p_\s,T_\s) u^\alpha_\welledge.  
$$
We first assume that both phases are present which implies that $T_{\rm sat} = (p_{\rm sat})^{-1}(p_\s)$ and that the liquid molar fraction is given by   
$$
c^\l_\s = { h^\g(p_\s,T_{\rm sat}) - {Q_{\s,e}^\omega \over Q_{\s,\mass}^\omega}  \over h^\g(p_\s,T_{\rm sat}) - h^\l(p_\s,T_{\rm sat})}, 
$$
from which we can compute $c^\g_\s = 1-c^\l_\s$ and the superficial velocities $u^\l_\welledge$ and $u^\g_\welledge$. 
The following alternatives are checked: 
\begin{itemize}
\item[] {\it Two-phase state}: if $0 < c^\l_\s <1$, the two-phase state is confirmed.
Using the slip law $u^\g = s^\g U_{\rm d}(s^\g) + s^\g C^0(s^\g) (u^g + u^\l)$, then 
$T_\s = T_{\rm sat}$ and the gas saturation $s^\g_\s$ is solution of the equation
$$
s^\g_\s U_{\rm d}(s^\g_\s) + s^\g_\s C^0(s^\g_\s) (u^g_\welledge + u^\l_\welledge) - u^\g_\welledge = 0. 
$$
In the no slip case, corresponding to $U_{\rm d} = 0$ and $C^0 = 1$, it reduces to $s_\s^\g = {u^\g_\welledge  \over u^g_\welledge + u^\l_\welledge}$. 
\item[] {\it Liquid state}: if $c^\l_\s \geq 1$, then only the liquid phase is present, we set $s_{\s}^{\l} = 1$, $s_{\s}^{\g} = 0$, and $T_\s$ is the solution of
  $$
  h^\l(p_\s,T_\s) = {Q_{\s,e}^\omega \over Q_{\s,\mass}^\omega}.
  $$
\item[] {\it Gas state}: if $c^\l_\s \leq 0$, then only the gas phase is present, we set $s_{\s}^{\l} = 0$, $s_{\s}^{\g} = 1$, and $T_\s$ is the solution of 
  $$
  h^\g(p_\s,T_\s) = {Q_{\s,e}^\omega \over Q_{\s,\mass}^\omega}.
  $$
\end{itemize}
We note that the same computations are done at the head node $\s_\omega$ based on the equations
$$
\begin{aligned}
 & \dsp  \sum_{\alpha\in \mathcal{P}}  |S^\omega_{\s_\omega}| \zeta^\alpha(p_{\s_\omega},T_{\s_\omega}) u^\alpha_\omega =  \sum_{\s''\in {\cal V}_\omega } q_{\s'',\mass}^{r\rightarrow \omega} := Q_{{\s_\omega},\mass}^\omega,   \\
&  \dsp \sum_{\alpha\in \mathcal{P}  }   |S^\omega_{\s_\omega}| h^\alpha(p_{\s_\omega},T_{\s_\omega}) \zeta^\alpha(p_{\s_\omega},T_{\s_\omega}) u^\alpha_\omega  =  \sum_{\s''\in {\cal V}_\omega}  q^{r\rightarrow \omega}_{\s'',e} := Q_{{\s_\omega},e}^\omega.   
\end{aligned}
$$
Then, the pressures $p_\s$, temperatures $T_\s$ and saturations $s_\s^\alpha$, $\s \in \mathcal{V}_\omega$ are used to compute the edge mean densities $\bar\rho^{\rm m}_\welledge$ defined by \eqref{eq_mum_rhom}. These mean densities are frozen in the momentum equation for the computation of the next time step. Together with the zero wall friction assumption (i) (which could be relaxed by an explicit approximation), it results that the next time step well pressures depend only on the well head node pressure $p_\omega$. This is the only well unknown implicitly coupled to the reservoir system combined with the single well equation reducing to the complementary constraints on the pair $\(\bar q_\omega - q_\omega, p_\omega-\bar p_\omega\)$  (first equation of \eqref{eq_monitoring}). Let us refer to see \cite{Xing2018} for more details about this type of single implicit unknown well model.

%%-----------------------------
%%      your bibliography
%%-----------------------------

\bibliographystyle{abbrv}
\bibliography{wellmodel-crossflow-Comp-Coats-v3}

\end{document}